\documentclass[final]{siamltex}
\usepackage{graphicx,psfrag}

\title{Influence of Intracellular Delay on the Dynamics of Hepatitis C Virus}

\author{Sandip Banerjee\thanks{Department of Mathematics, Indian Institute of Technology Roorkee,
Roorkee, Uttarakhand, 247667, India ({\tt sandofma@iitr.ac.in}).}
        \and Ram Keval\thanks{Department of Mathematics, Indian Institute of Technology
        Roorkee, Roorkee, Uttarakhand, 247667, India ({\tt ramkeval@gmail.com}).} \and Sunita Gakkhar\thanks{Department of Mathematics,
        Indian Institute of Technology Roorkee, Roorkee, Uttarakhand, 247667, India ({\tt sungkfma@iitr.ac.in}).}}

\begin{document}

\maketitle

\begin{abstract}
In this paper we present a delay induced model for hepatitis C virus incorporating the healthy and infected hepatocytes as well as infectious and noninfectious virions. The model is mathematically analyzed and characterized, both for the steady states and the dynamical behavior of the model. It is shown that time delay does not affect the local asymptotic stability of the uninfected steady state. However, it can destabilize the endemic equilibrium, leading to Hopf bifurcation to periodic solutions with realistic data sets. The model is also validated using 12 patient data obtained from the study, conducted at the University of Sao Paulo Hospital das clinicas.
\end{abstract}



\pagestyle{myheadings}
\thispagestyle{plain}
\markboth{Sandip Banerjee, Ram Keval and Sunita Gakkhar}
{Influence of Intracellular Delay on the Dynamics of Hepatitis C Virus}

\section{Introduction}

An overwhelming number of people globally are infected with the hepatitis C virus (HCV). Chronic HCV often requires
treatment in the form of liver transplant \cite{Dixit08}. Unfortunately, most of the cases of chronic HCV lead to fatalities. The traditional treatment with pegylated interferon (IFN) has met with moderate success. The current front line therapy for HCV involves a combination of IFN and ribavirin \cite{Dixit08, Feld05}. Even then this combination therapy elicits success in about half of the treated cases \cite{Dixit08} and is quite expensive. It is therefore imperative to study the mechanism of viral
dynamics and the precise extent of therapeutic response. The eventual goal would be to improve upon the therapeutic
protocols \cite{Feld05} keeping in mind the goal of minimizing the viral load while also minimizing the toxic side
effects \cite{ChakrabartyJBS}. As already pointed out the therapeutic effectiveness of IFN is moderate but patients show
remarkable progress when it is combined with ribavirin \cite{Feld05, Graci06}. The response to the combination therapy
has been varied in nature \cite{Dixit08, Feld05}. In some cases, there has been a sustained fall of the virus below
detection levels, whereas in other cases there has either been a fall in the viral load followed by a relapse after
cession of treatment as well as cases of complete non response \cite{Dixit08, Feld05}. Even in successful cases, the
viral load has exhibited biphasic and triphasic levels of decline \cite{Dahari07}.

Mathematical models which are simple but elegant have played a pivotal role (especially in recent years)
in the understanding of the mechanism of viral dynamics in the cases of influenza, hepatitis B and
human immunodeficiency virus (HIV) \cite{Dixit05, Nowak96, Perelson99, Baccam06, Lewin01}. One of the earliest realistic model to reflect
the viral dynamics of HCV was a work by Neumann et al. that appeared in 1998 \cite{Neumann98}.
The model \cite{Neumann98} comprised of a simple but coupled system of three ordinary differential equations (ODEs).
It involved the mechanism of the uninfected and the productively infected target cells (which in case of HCV are the
hepatocytes). The treatment involved only IFN and did not include ribavirin. The hepatocytes were assumed to have a
constant source of supply and undergo a natural death. Some of these hepatocytes were rendered into an infected state
by the HCV resulting in infected hepatocytes, some of which were lost due to immune response. The virions were following
a birth (due to infected hepatocytes) and a natural death process \cite{Dixit08}. The authors then incorporated into the
model the effectiveness of IFN on the viral load. It was observed that the IFN was primarily playing the role of
rendering the virions ineffective, while having a minimal role in controlling the transition of hepatocytes from an
uninfected to an infected state. An assumption in the contrary to the previous statement, resulted in the model showing
results of a single phase delay which contradicted the clinical results. However, when the assumption of IFN being mainly
used to target the virus, was used in the model, the results showed a close matching with the experimental
results \cite{Graci06, Neumann98} of biphasic decline in the viral load. Once this model was validated for
\textit{in vivo} HCV dynamics it was used to estimate key parameter values of the model. The drawback of this model
was its inability to capture the long-term behavior of dynamics as well as it's failure to take into account
the consequences of usage of ribavirin.

The motivation of a subsequent work was to address the issues that were left open for deliberations
in \cite{Neumann98}. Dixit et al. \cite{Dixit04}, in an endeavor to capture the effects of ribavirin
in the dynamics, presented an improved model. In this model the presence of ribavirin created two categories
of viruses (one infectious and the other noninfectious) owing to the understanding that the predominant role
of ribavirin was to render a fraction of the virus into a noninfectious state. The most interesting finding of this
work was the observation in numerical simulations that ribavirin has a little role to play in the first phase of a
biphasic decline of viral load. It also affects the second phase and that too only when the effectiveness of IFN is
moderate. The biological explanation of this as provided by the authors \cite{Dixit04} is that when IFN is effective,
it can control the virions and hence ribavirin is left with a little role to play. It is when IFN is minimally effective
that ribavirin reduces the viral load not by arresting the virus which in turn would have affected the production of
infected hepatocytes, but because of the fact that viruses have been rendered noninfectious \cite{Dixit08}.
 This model was also successful in capturing some of the long-term effects.

While both these models \cite{Neumann98, Dixit04} are widely regarded as seminal work in mathematical modeling of HCV dynamics, they did not fully address the clinical observations of triphasic decline and more importantly, the critical issue of non-response (which unfortunately happens in a majority of the chronic cases). These limitations of \cite{Neumann98, Dixit04} were addressed in a recent work by Dahari et al. \cite{Dahari07}. This model was an extension of the previous models \cite{Neumann98, Dixit04}. The model for the first time included the mechanism of the liver by including a logistic growth term for the hepatocytes (both uninfected and infected) and a maximum carrying capacity which put a limit on the proliferation of hepatocytes. They were able to obtain an elegant, yet simple critical effectiveness threshold beyond which one cannot attain sustained long time viral decay. This threshold which was a  function of both viral and host parameters, was able to explain the rationale being the non-effectiveness of the  treatment. This homeostatic cell proliferation terms were also able to address mathematically the causes behind the triphasic decline observed in some patients. This model, however, still has aspects that have to be addressed in terms of clinical validation, since the estimation of parameters from liver mechanism still remains a gray area. A similar model has been proposed by Banerjee et al. \cite{Banerjee13}, where dynamical behavior of the HCV due to separate efficacies of interferon and ribavirin hold the key factor.

Time delays in connection with viral dynamics, namely, HIV, Lymphocytic  choriomeningitis virus (LCMV), Hepatitis B are always an interesting features of study among the researches as models of viral infection dynamics with delays may have complicated impact on the dynamics of a system. Also, a delay in the viral production of infected cells may drastically change the estimates of viral clearance \cite{Nelson02}. Delay may cause the loss of stability and can bifurcate various periodic solutions \cite{Gourley08}. Clearly in viral dynamics, the time for viral infection is not instantaneous. Development of the mathematical model and the associated analytical and numerical techniques that incorporate discrete time delay into the model of hepatitis C viral dynamics is the primary motivation here.

Section 2 gives the formulation of the model in details. Section 3 deals with the linear stability analysis and bifurcation analysis, including the estimation of the length of delay to preserve stability, direction and stability of Hopf bifurcation. Section 4 gives the numerical results and its biological implications and the paper ends with a discussion.

\section{Model Formulation}

The proposed model is the amalgamation and modification of the models that appear in \cite{ChakrabartyJBS, Banerjee13}, which is represented by the system of coupled ODEs:
\begin{eqnarray}\label{eq1}
\frac{dT}{dt}& = & s + r T \left(1-\frac{T+I}{T_{\max}}\right)-d_{1}T-(1-c\eta_{1})\alpha V_{I}T, \nonumber\\
\frac{dI}{dt}& = &(1-c\eta_{1})\alpha V_{I}T -d_{2} I,\nonumber\\
\frac{dV_{I}}{dt}& = & \left(1-\frac{\eta_{r}+\eta_{1}}{2}\right)\beta I- d_{3}V_{I},\\
\frac{dV_{NI}}{dt}& = & \left(\frac{\eta_{r}+\eta_{1}}{2}\right)\beta I -d_{3}V_{NI}.\nonumber
\end{eqnarray}
Here, $T$ denotes the target cells (uninfected hepatocytes), $I$ are infected hepatocytes (productively infected cells), $V_{I}$ and $V_{NI}$ represent the infectious and non-infectious viruses (virions) respectively. The first equation projects the growth of uninfected hepatocytes from a constant source at a rate s, followed by an intrinsic logistic growth term at a rate r. Please note that the logistic part is of the form $r T (1-(T+I)/T_{\max})$, where $T_{\max}$ is the carrying capacity of the hepatocytes (both infected as well as uninfected). This is because of the fact that as the threshold value $T_{\max}$ is attained, the infected and the non-infected hepatocytes are not differentiated but considered as a whole. As already pointed out in \cite{Dahari07}, the logistic term is incorporated to capture the homeostatic cell proliferation. There are two drugs acting on the system, namely, interferon and ribavirin. At this point we note that $\eta_1$ is the efficacy of interferon in blocking the release of new virions. The term $(1-c \eta_{1})$ ($0 < c < 1$) reflects the ineffectiveness of interferon (which is much smaller than its effectiveness in blocking the virions and which is why we take $(0 < c < 1$)) in inhibiting hepatocytes conversion from the uninfected to the infected one.

The second equation gives the dynamics of infectious hepatocytes, where the uninfected hepatocytes (T) are infected by the infectious virions ($V_I$) at the rate $\alpha$ and the whole term $(1-cn_{1})\alpha V_{I} T$ shifts from the uninfected class to the infected one; ($-d_2 I$) is the natural death of the infected hepatocytes (I).

The combined effect of the interferon and ribavirin has been incorporated in the last two equations. Due to the combined efficacy, the infectious hepatocytes replicates into infectious ($V_I$) as well as non-infectious ($V_{NI}$) virions. The term $(1-(\eta_{r} + \eta_{1})/2)\beta I$ determines the growth of infectious virions due to non-efficacy of the combined drugs: interferon and ribavirin, due to which the infected cells are replicating virions at a rate $\beta$. But, due the efficacy of the combined effect of interferon and ribavirin, the infected cells are replicating non-infectious virions at a rate $\beta$, which is given by the term $\left((\eta_{r} + \eta_{1})/2\right)\beta I$, which also points out that the combined effect of the drugs interferon and ribavirin are successful in rendering the virions from an infected state to an uninfected one. Finally, the infectious and noninfectious viruses are assumed to have the same natural death rate of $d_{3}$. Note that here $0 < \eta_{1} < 1$ and $0< \eta_{r} < 1$ lead to $0<\left((\eta_{r} + \eta_{1})/2\right)<1$.

It has been almost invariably observed that in viral dynamics, there is a time delay (say, $\tau$) in conversion of the target cells (the CD4 T-cells in case of HIV \cite{Nelson02} and hepatitis B viral dynamics \cite{Gourley08} for instance) from an uninfected to an infected stage. This delay in transformation induces delay in the growth of infected cells. Hence, the model, which is represented by a system of coupled ODEs (equation(\ref{eq1})) is converted to a system of delay differential equations (DDEs):

\begin{eqnarray}\label{eq2}
\frac{dT}{dt}& = & s+r T \left(1-\frac{T+I}{T_{\max}}\right)-d_{1}T -(1-c\eta_{1})\alpha V_{I} T,\nonumber\\
\frac{dI}{dt}& = &(1-c\eta_{1})\alpha V_{I}(t-\tau)T(t-\tau) -d_{2} I,\nonumber\\
\frac{dV_{I}}{dt}& = &\left(1-\frac{\eta_{r}+\eta_{1}}{2}\right)\beta I- d_{3}V_{I},\\
\frac{dV_{NI}}{dt}& = &\left(\frac{\eta_{r}+\eta_{1}}{2}\right)\beta I -d_{3}V_{NI}.\nonumber
\end{eqnarray}
The system is subjected to initial conditions:
%\begin{equation}
%T(\theta)=\psi_1(\theta),~I(\theta) = \psi_2(\theta),~V_I(\theta) = \psi_3(\theta),~ V_{NI}(\theta) = \psi_4(\theta),~\theta \in [-\tau,~0],~\psi_i(0) \geq 0,~i=1,2,3,4
%\end{equation}
\begin{eqnarray}\label{inicialcond}
\nonumber ~~~~~~~~T(\theta)&=&\psi_1(\theta),~I(\theta) = \psi_2(\theta),~V_I(\theta) = \psi_3(\theta),~ V_{NI}(\theta) = \psi_4(\theta),~\theta \in [-\tau,~0],\\
\psi_i(0) & \geq & 0,~i=1,2,3,4
\end{eqnarray}

where $\left(\psi_1(\theta), \psi_2(\theta), \psi_3(\theta), \psi_4(\theta)\right) \in C\left([-\tau,~0],R^{4}_{+0}\right)$, the Banach space of continuous functions, mapping the interval $[-\tau,~0]$ into $R^{4}_{+0}$. We define $R^{4}_{+0}$ and $R^4_{+}$ (the interior of $R^4_{+0}$) as:
\begin{eqnarray}
\nonumber R^4_{+0} &=& \left((T,I,V_I,V_{NI}):T,I,V_I,V_{NI}\geq0\right)\\
 \nonumber R^4_{+} &=& \left((T,I,V_I,V_{NI}):T,I,V_I,V_{NI}>0\right)
\end{eqnarray}
This model of DDEs will form the fundamental mathematical framework for this paper.

\section{Qualitative analysis of the model}
\subsection{Equilibria and linear stability analysis}

There exit two positive equilibrium points of (\ref{eq1}), namely,
\begin{enumerate}
\item $E_{1} : $
\[\widehat{T}= \frac{T_{\max}}{2 r}\left[(r-d_{1})+\sqrt{(r-d_{1})^{2}+\frac{4 r s}{T_{\max}}}~\right],~~r > d_{1},~ I=0,~ V_{I}=0,~ V_{NI}=0.\]

\item $E_{2} : $
\begin{eqnarray}
\nonumber T^{*}&=&\frac{d_{2}d_{3}}{{(1-c\eta_{1})\left(1-\frac{\eta_{r}+\eta_{1}}{2}\right)\alpha \beta}}, ~~~~I^{*}=\left(\frac{s R_{0} T_{\max}+r \widehat{T}^{2}}{r \widehat{T}+d_{2}R_{0}T_{\max}}\right)\left(1-\frac{1}{R_{0}}\right),\\
\nonumber V_{I}^{*}&=&\left(1-\frac{\eta_{r}+\eta_{1}}{2}\right)\frac{\beta I^{*}}{d_{3}}, ~~~~~~~~V_{NI}^{*}=\left(\frac{\eta_{r}+\eta_{1}}{2}\right)\frac{\beta I^{*}}{d_{3}}.
\end{eqnarray}
\end{enumerate}
Here $R_{0}=\frac{\widehat{T}}{T^{*}}$ is the basic reproduction number. $E_{2}$ exists if $R_{0}>1.$ We assume that $s \leq d_{1} T_{\max}$, so that the model is physiologically realistic.

The characteristic equation about any equilibrium point $\bar{E}(\bar{T}, \bar{I}, \bar{V}_{I}, \bar{V}_{NI})$ is given by
\begin{equation}\label{charact}
    \left|
  \begin{array}{cccc}
    r-d_{1}- \frac{2 r \bar{T}}{T_{\max}}-\frac{r \bar{I}}{T_{\max}}-(1-c\eta_{1})\alpha \bar{V}_{I}-\lambda& -\frac{r \bar{T}}{T_{\max}} & -(1-c\eta_{1})\alpha \bar{T} & 0 \\
    (1-c\eta_{1})\alpha \bar{V}_{I} e^{-\lambda \tau} & -d_{2}-\lambda&(1-c\eta_{1})\alpha \bar{T} e^{-\lambda \tau}  & 0 \\
    0 & \left(1-\frac{\eta_{r}+\eta_{1}}{2}\right)\beta & -d_{3}-\lambda &0 \\
    0 & \left(\frac{\eta_{r}+\eta_{1}}{2}\right)\beta & 0 & -d_{3}-\lambda
  \end{array}
\right|=0.
\end{equation}

For the equilibrium point $E_{1}=(\widehat{T}, 0, 0, 0)$, (\ref{charact}) reduces to
\begin{eqnarray}\label{infcharac}
\nonumber && \left(\lambda + \frac{2 r}{T_{\max}}\sqrt{(r-d_{1})^{2}+\frac{4 r s}{T_{\max}}}~\right)(\lambda+d_{3})
\left[\lambda^{2}+ (d_{2}+d_{3})\lambda+d_{2}d_{3}\left(1-R_0 e^{-\lambda \tau}\right)\right]=0.
\end{eqnarray}
Clearly, for $\tau = 0,$ $E_{1}$ is locally asymptotically stable for $R_{0}<1$ and unstable for $R_{0}>1$. For $\tau > 0,$ two eigenvalues are negative and the other two are given by
\begin{equation}\label{infcharac1}
    \lambda^{2}+ (d_{2}+d_{3})\lambda+d_{2}d_{3}\left(1-R_0 e^{-\lambda \tau}\right)=0.
\end{equation}

\textbf{Case I ($\tau > 0, R_0 > 1$)}:~ Let $F_1(\lambda) = \lambda^{2}+ (d_{2}+d_{3})\lambda$ and $F_2(\lambda) = d_{2}d_{3}\left(R_0 e^{-\lambda \tau}-1\right)$, so that $F_1(\lambda) = F_2(\lambda)$. Clearly, $F_1(0) = 0$ and $F_1(\lambda) \rightarrow\infty$ as $\lambda \rightarrow \infty$; whereas $F_2(\lambda)$ is decreasing in $\lambda$. Furthermore, $F_2(0) = d_{2}d_{3}\left(R_0 -1\right) > 0$ (since, $R_0 > 1$). Therefore, for some positive $\lambda$, the two functions $F_1(\lambda)$ and $F_2(\lambda)$ must intersect, implying the fact that the equation (\ref{infcharac1}) has a positive root and we conclude that $E_1$, the disease free steady state is unstable for positive delay and $R_0 > 1$. \\

\textbf{Case II ($\tau > 0, R_0 < 1$)}:~ Our expectation in this case is (\ref{infcharac1}) possesses negative roots and the disease free equilibrium is stable. Clearly,  $\lambda^{2}+ (d_{2}+d_{3})\lambda = d_{2}d_{3}\left(R_0 e^{-\lambda \tau}-1\right)$ does not possess a non-negative real root because $F_1(\lambda)$ is an increasing function for $\lambda \geq 0$, whereas $F_2(\lambda)$ is a decreasing function with $F_2(0) = R_0 - 1 < 0$. Therefore, we conclude that (\ref{infcharac1}) will have roots with non-negative real parts if the roots are complex and should have been obtained from a pair of complex conjugate roots, which crosses the imaginary axis \cite{Culshaw00}. This implies that for some $\tau \geq 0$, (\ref{infcharac1}) must have a pair of purely imaginary roots. Without any loss of generality, let $\lambda = i \omega$ ($\omega > 0$) be the root of (\ref{infcharac1}). Substituting in (\ref{infcharac1}) and separating the real and imaginary parts we get,
\begin{eqnarray}
  \nonumber -\omega^{2}+d_{2} d_{3} &=& R_0 d_2 d_3 \cos (\omega \tau) ,\\
  \nonumber (d_{2}+ d_{3})\omega &=& R_0 d_2 d_3  \sin (\omega \tau).
\end{eqnarray}
Eliminating the trigonometric function we get,
\begin{eqnarray} \label{eqomega}
(\omega^{2})^2 + (d_{2}^{2}+ d_{3}^{2}) \omega^{2} + d_{2}^{2} d_{3}^{2} \left(1-R_0^{2}\right) = 0.
\end{eqnarray}
Clearly, in (\ref{eqomega}), the sum of the roots are negative and the product of the roots are positive (as $R_0 < 1$) and consequently, (\ref{eqomega}) does not have any positive real roots, implying that there does not exist any $\omega$ for which $i \omega$ is a solution of (\ref{infcharac1}). Following Rouch\'{e}'s theorem \cite{Dieudonne60} (Theorem 9.17.4), we conclude that all the eigen values of the characteristic equation (\ref{infcharac1}) have a negative real part for all $\tau \geq 0$. Therefore, $E_1$ is locally asymptotically stable for $R_0 < 1$.

\subsection{\textbf{Critical drug efficacy}}
The critical drug efficacy for this model is defined as (for details see \cite{Banerjee13})
\begin{eqnarray*}
\eta_c = 1 - \frac{2 r d_{2}d_{3}}{\alpha \beta [(r-d_1) k + \sqrt{(r-d_1)^2 k^2+ 4 r s k}]} = 1 - \frac{T_0^{*}}{\widehat{T}}
\end{eqnarray*}
where $T_{0}^{*} = \frac{d_2 d_3}{\alpha \beta}$ is the number of uninfected hepatocytes in an infected person before treatment (obtained by putting $\eta_1 = \eta_r = 0$ in $T^*$) and $\widehat{T}$ is the steady state of hepatocytes in an uninfected individual. $\eta_c$ gives a measure for the threshold of the efficacies of drugs, interferon and ribavirin. If the combined efficacies exceed the critical value, HCV are eradicated whereas if $\eta < \eta_c$, HCV reaches a new steady state, lower than the previous steady state value \cite{Dahari07, Banerjee13}.

\subsection{Hopf Bifurcation Analysis}

For the equilibrium point $E_{2}=(T^*,I^*,V_I^*, V_{NI}^*)$, (\ref{charact}) reduces to

\begin{equation}\label{delaycharacteristiceq1}
(\lambda +d_3)[\lambda^3 + a_0 \lambda^2 + a_1 \lambda + a_2 + ( b_1 \lambda + b_2) e^{-\lambda \tau}]=0
\end{equation}
where,
\begin{eqnarray}
 \nonumber a_{0}&=& d_{2}+d_{3}+\frac{s}{T^{*}}+\frac{r T^{*}}{T_{\max}}, ~a_{1}=d_{2}d_{3}+(d_{2}+d_{3})\left(\frac{s}{T^{*}}+\frac{r T^{*}}{T_{\max}}\right),\\
\nonumber a_{2}&=& d_{2}d_{3}\left(\frac{s}{T^{*}}+\frac{r T^{*}}{T_{\max}}\right),\\
\nonumber b_{1}&=& \frac{d_{2}r I^{*}}{T_{\max}}-(1-c\eta_{1}) \left(1-\frac{\eta_{r}+\eta_{1}}{2}\right)\alpha \beta T^{*}~~=\frac{d_{2} r I^{*}}{T_{\max}}-d_{2} d_{3},\\
 \nonumber b_2 &=& \frac{d_{2} d_{3} r I^{*}}{T_{\max}}+\left(d_{2} I^{*}-s-\frac{r T^{*2}}{T_{\max}}\right)(1-c\eta_{1}) \left(1-\frac{\eta_{r}+\eta_{1}}{2}\right)\alpha \beta\\
\nonumber &=& d_{2} d_{3} I^{*}\left(\frac{r}{T_{\max}}+\frac{d_{2}}{T^{*}}\right)-d_{2}d_{3}\left(\frac{s}{T^{*}}+\frac{r T^{*}}{T_{\max}}\right).
\end{eqnarray}

Clearly, the characteristic equation (\ref{delaycharacteristiceq1}) has one negative eigenvalue, namely, $-d_3$ and we are left with
\begin{equation}\label{delaycharacteristiceq}
\lambda^3 + a_0 \lambda^2 + a_1 \lambda + a_2 + ( b_1 \lambda + b_2) e^{-\lambda \tau}=0,
\end{equation}
a transcendental equation (\ref{delaycharacteristiceq}) which has infinitely many eigenvalues and for which the classical Routh-Hurwitz condition for stability analysis fails.\\

In the absence of delay ($\tau=0$), the characteristic equation (\ref{delaycharacteristiceq}) becomes
\begin{equation}\label{characteristiceq}
\lambda^3 + a_0 \lambda^2 + (a_1 + b_1) \lambda + (a_2 + b_2)=0,
\end{equation}
By Routh-Hurwitz criteria, (\ref{characteristiceq}) will have negative real parts if $a_0 >0$, $a_{1} + b_{1}>0$, $a_{2} + b_{2}>0$ and $(a_{1}+b_{1})a_{0}-(a_{2} + b_{2}) > 0$. Clearly, the first three conditions hold. Therefore, the equilibrium point $E_{2}$ is asymptotically stable in the absence of delay ($\tau =  0$), provided
\begin{eqnarray}
\label{condition1}
\left[(d_{2}+d_{3})\left(\frac{s}{T^{*}}+\frac{r T^{*}}{T_{\max}}\right)+\frac{d_{2} r I^{*}}{T_{\max}}\right] \times \left(d_{2}+d_{3}+\frac{s}{T^{*}}+\frac{r T^{*}}{T_{\max}}\right)
\\\nonumber-d_{2} d_{3} I^{*}\left(\frac{r}{T_{\max}}+\frac{d_{2}}{T^{*}}\right)>0
\end{eqnarray}

We now want to investigate the existence of purely imaginary roots of (\ref{delaycharacteristiceq}). Putting $\lambda = i \omega$ ($\omega$ is taken as positive) in the equation (\ref{delaycharacteristiceq}) and solving for real and imaginary parts, we have the system of transcendental equations as
\begin{eqnarray}\label{transcendentaleq1}
a_{0}\omega^2 -a_{2} &=& b_{2} \cos(\omega \tau)+ b_{1}\omega \sin(\omega \tau)\\
\label{transcendentaleq2}
\omega^{3}- a_{1} \omega &=& b_{1}\omega  \cos(\omega
\tau)-b_2  \sin(\omega \tau)
\end{eqnarray}
Squaring and adding (\ref{transcendentaleq1}) and (\ref{transcendentaleq2}), we obtain,
\begin{equation}\label{transcendentaleq3}
\omega^6 + A_{1}\omega^4 +A_{2}\omega^2+A_{3}=0,
\end{equation}
with
\begin{equation}\label{transcendentaleq4}
A_{1}=(a_0^2-2 a_1),~~~~ A_{2}= (a_1^2-b_1^2 - 2 a_0 a_2),~~~~ A_{3} = a_2^2-b_2^2.
\end{equation}
The simplest assumption that (\ref{transcendentaleq3}) will have a positive root is $A_1 > 0$ and $A_3 < 0$. Clearly, $A_{1}=d_{2}^{2}+d_{3}^{2}+\left(\frac{s}{T^{*}}+\frac{r T^{*}}{T_{\max}}\right)^{2}$ is positive. Therefore, for a positive root of the equation (\ref{transcendentaleq3}), we simply conclude that $a_2-b_2 < 0$, that is,
\begin{eqnarray*}
2 d_{2}d_{3}\left(\frac{s}{T^{*}}+\frac{r T^{*}}{T_{\max}}\right)-d_{2} d_{3} I^{*}\left(\frac{r}{T_{\max}}+\frac{d_{2}}{T^{*}}\right) < 0
\end{eqnarray*}
Therefore, we can say that there exists a unique positive root $\omega_0$ satisfying (\ref{transcendentaleq3}), that is, the characteristic equation (\ref{delaycharacteristiceq1}) will have purely imaginary roots of the form $\pm i \omega$. From (\ref{transcendentaleq1}) and (\ref{transcendentaleq2}), solving for $\cos(\omega \tau)$,  we get
\begin{eqnarray}
\nonumber \cos(\omega \tau)=\frac{(a_{0} \omega^{2}-a_2)b_2 + (\omega^{3}-a_1 \omega) b_1 \omega}{b_2^{2}+ (b_1 \omega)^2}.
\end{eqnarray}
Then $\tau_j$ corresponding to $\omega_{0}$ is given by
\begin{eqnarray}\label{proposition1eq4}
~~~~~\tau_j=\frac{1}{\omega_0}\arccos\left[\frac{(a_{0} \omega_0^{2}-a_2)b_2 + (\omega_0^{3}-a_1 \omega_0) b_1 \omega_0}{b_2^{2}+ (b_1 \omega_0)^2} \right]+\frac{2j \pi}{\omega_0},~~~ j = 0,1,2,.....
\end{eqnarray}
For $\tau = 0, E_{2}$ is locally asymptotically stable, provided (\ref{condition1}) holds, therefore, by Butler's Lemma \cite{Freedman83}, $E_{2}$ remains locally asymptotically stable for $\tau <\tau_0$.\\

To establish the Hopf bifurcation at $\tau = \tau_0$, we have to show the transversality condition
\begin{eqnarray}
\nonumber \left[\frac{d(\textrm{Re}\lambda)}{d\tau}\right]_{\tau~=~\tau_j}>0.
\end{eqnarray}

Now, differentiating (\ref{delaycharacteristiceq}) with respect to $\tau$, we get
\begin{eqnarray}
\nonumber \left[(3\lambda^2+2a_0\lambda+a_1)+e^{-\lambda \tau} b_1-\tau e^{-\lambda \tau}(b_1\lambda+b_2)\right]\frac{d\lambda}{d\tau} =(b_1\lambda+b_2)e^{-\lambda \tau}\lambda.
\end{eqnarray}
This implies
\begin{eqnarray}
\nonumber \left[\frac{d\lambda}{d\tau}\right]^{-1}&=&\frac{3\lambda^2+2a_0\lambda+a_1}{\lambda e^{-\lambda \tau} (b_1\lambda+b_2)}+\frac{b_1}{\lambda (b_1\lambda+b_2)} -\frac{\tau}{\lambda}\\
\nonumber &=& \frac{3\lambda^2+2a_0\lambda+a_1}{-\lambda (\lambda^3+a_0 \lambda^2+a_1\lambda+a_2)}+\frac{b_1}{\lambda (b_1\lambda+b_2)}-\frac{\tau}{\lambda}\\
\nonumber &=&\frac{2\lambda^3+a_0\lambda^2-a_2}{-\lambda^2 (\lambda^3+a_0 \lambda^2+a_1\lambda+a_2)}+\frac{- b_2}{\lambda^2 (b_1\lambda+b_2)}-\frac{\tau}{\lambda}
\end{eqnarray}
Therefore,
\begin{eqnarray}
\nonumber \Theta&=&\textmd{sign} \left[Re \left(\frac{2\lambda^3+a_0\lambda^2-a_2}{-\lambda^2 (\lambda^3+a_0 \lambda^2+a_1\lambda+a_2)}+\frac{- b_2}{\lambda^2 (b_1\lambda+b_2)}-\frac{\tau}{\lambda}\right)\right]_{\lambda ~=~ i \omega_0}\\
\nonumber &=&\frac{1}{\omega_0^2}\textmd{sign} \left[Re \left(\frac{(a_2+a_0 \omega_0^2)+i 2\omega_0^3}{(a_0\omega_0^2-a_2)+i(\omega_0^3-a_1\omega_0)} +\frac{b_2}{b_2+i (b_1\omega_0)}\right)\right]\\
\nonumber &=&\frac{1}{\omega_0^2} \textmd{sign} \left[\frac{(a_2+a_0\omega_0^2)(a_0\omega_0^2-a_2) +2\omega_0^3(\omega_0^3-a_1\omega_0)}{(a_0\omega_0^2-a_2)^2 +(\omega_0^3-a_1\omega_0)^2} +\frac{b_2^{2}}{b_2^2+ (b_1\omega_0)^2}\right]\\
\nonumber &=&\frac{1}{\omega_0^2}\textmd{sign} \left[\frac{2\omega_0^6+(a_0^2-2a_1) \omega_0^4+(b_2^2-a_2^2)}{b_2^2+(b_1\omega_0)^2}\right]
\end{eqnarray}
Clearly,  $A_{1} = a_0^2-2a_1$ is positive and by the assumption $A_{3}=a_2^2-b_2^2$ is negative, implying
\begin{eqnarray}
\nonumber \left[\frac{d(Re \lambda)}{d\tau}\right]_{\omega ~=~ \omega_0,~ \tau ~=~ \tau_j}>0.
\end{eqnarray}
We summarize the preceding analysis in the form of a theorem.

\begin{theorem}\label{dth5}
Suppose\\
(i) $R_0 > 1$.\\
(ii) $\left[(d_{2}+d_{3})\left(\frac{s}{T^{*}}+\frac{r T^{*}}{T_{\max}}\right)+\frac{d_{2} r I^{*}}{T_{\max}}\right] \left(d_{2}+d_{3}+\frac{s}{T^{*}}+\frac{r T^{*}}{T_{\max}}\right)-d_{2} d_{3} I^{*}\left(\frac{r}{T_{\max}}+\frac{d_{2}}{T^{*}}\right)>0$ is satisfied and\\
(iii) the largest positive simple root of (\ref{transcendentaleq3}) is $\omega_0$,\\
then $E_{2}$, the endemic equilibrium point of the delay model is asymptotically stable when $\tau<\tau_0$ and unstable when $\tau>\tau_0$, where,
\begin{eqnarray*}
~~~~~\tau_0=\frac{1}{\omega_0}\arccos\left[\frac{(a_{0} \omega_0^{2}-a_2)b_2 + (\omega_0^{3}-a_1 \omega_0) b_1 \omega_0}{b_2^{2}+ (b_1 \omega_0)^2} \right].
\end{eqnarray*}
At $\tau$ passes through the critical point $\tau=\tau_0$, a family of periodic solutions bifurcates from $E_{2}$, that is, a Hopf bifurcation occurs at $\tau=\tau_0$ .
\end{theorem}

\subsection{Estimation of the length of delay to preserve stability}

It should be noted that the stability of bifurcating periodic orbits cannot be determined from the previous analysis, that is, in the neighborhood of $\tau_0$, existence of periodic solutions will depend on either for $\tau > \tau_0$ or for $\tau < \tau_0$. While investigating the stability of bifurcating periodic orbits, we now try to estimate the maximum length of delay to maintain the stability of the limit cycle. For this, we consider the space of real valued continuous functions, defined on $[-\tau, ~\infty)$, satisfying the initial conditions on $[-\tau,~ 0].$

Let $T(t)= T^{*}+X(t), I(t)= I^{*}+Y(t), V_{I}(t)= V_{I}^{*}+Z(t), V_{NI}(t)= V_{NI}^{*}+W(t)$. Linearizing the system (\ref{eq2}) about the equilibrium point $E_{2}(T^{*}, I^{*}, V_{I}^{*}, V_{NI}^{*})$, we get
\begin{eqnarray}\label{linsystem}
\nonumber \frac{dX}{dt}&=& \left(r-d_{1}-\frac{2 r T^{*}}{T_{\max}}-\frac{r I^{*}}{T_{\max}}-(1-c \eta_{1})\alpha V_{I}^{*}\right) X(t)-\frac{r T^{*}}{T_{\max}}Y(t)\\
\nonumber &-&(1-c \eta_{1})\alpha T^{*} Z(t)\\
 \nonumber \frac{dY}{dt} &=& (1-c \eta_{1})\alpha V_{I}^{*} X(t-\tau)-d_{2}Y(t)+(1-c \eta_{1})\alpha T^{*} Z(t-\tau)\\
\frac{dZ}{dt} &=& \left(1-\frac{\eta_{1}+\eta_{r}}{2}\right)\beta Y(t)-d_{3} Z(t)\\
\nonumber \frac{dW}{dt} &= &\left(\frac{\eta_{1}+\eta_{r}}{2}\right)\beta Y(t)-d_{3} W(t).
\end{eqnarray}
Using Laplace transform of the linearized system (\ref{linsystem}), we obtain,
$$
\left\{ \begin{array}{ll}
\left(\delta-r+d_{1}+\frac{2 r T^{*}}{T_{\max}}+\frac{r I^{*}}{T_{\max}}+(1-c \eta_{1})\alpha V_{I}^{*}\right)L_{X}(\delta)=-\frac{r T^{*}}{T_{\max}} L_{Y} (\delta)\\
-(1-c \eta_{1})\alpha T^{*} L_{Z}(\delta)+X(0)\\
(\delta+d_{2})L_{Y}(\delta) = (1-c \eta_{1})\alpha e^{-\delta \tau} V_{I}^{*}L_{X}(\delta)\\
+(1-c \eta_{1})\alpha e^{-\delta \tau}
V_{I}^{*}K_{1}(\delta)+(1-c \eta_{1})\alpha e^{-\delta \tau} T^{*}L_{Z}(\delta)+(1-c \eta_{1})\alpha e^{-\delta\tau}T^{*}K_{2}(\delta)+Y(0)\\
(\delta+d_{3})L_{Z}(\delta)=\left(1-\frac{\eta_{1}+\eta_{r}}{2}\right)\beta L_{Y}(\delta)+Z(0)\\
(\delta+d_{3})L_{W}(\delta)=\left(\frac{\eta_{1}+\eta_{r}}{2}\right)\beta L_{Y}(\delta)+W(0)
 \end{array} \right.
$$
where, $ K_{1}(\delta) = \int^{0}_{-\tau} e^{-\delta t}X(t) dt, K_{2}(\delta) = \int^{0}_{-\tau}e^{-\delta t} Z(t)dt$ and $L_{X}, L_{Y}, L_{Z}$ and $L_{W}$ are the Laplace transforms of $X(t), Y(t), Z(t)$ and $W(t)$ respectively.

Applying Nyquist criterion, we obtain the conditions for local asymptotic stability of $E_{2}(T^{*}, I^{*}, V_{I}^{*}, V_{NI}^{*})$ as \cite{Freedman83}
\begin{eqnarray}\label{eqn1}
\textrm{Im} H(i \mu_{0})> 0,\\
\label{eqn2}
\textrm{Re} H(i \mu_{0})= 0.
\end{eqnarray}
where $H(\delta)= (\delta^{3}+a_{0}\delta^{2}+a_{1}\delta+a_{2})+e^{-\delta \tau}(b_{1}\delta+b_{2})$, $\mu_{0}$ being the smallest positive root of (\ref{eqn2}). It has already been proved that in the absence of delay, $E_2$ is locally asymptotically stable, provided (\ref{condition1}) holds. Then, by virtue of Burlet's lemma \cite{Freedman83}, we can state that for sufficiently small $\tau > 0$, all the  eigenvalues will continue to have negative real parts (by continuity), provided there is a guarantee that no eigenvalues with positive real parts bifurcates from infinity as $\tau$ increases from zero.\\

Inequality (\ref{eqn1}) and equation (\ref{eqn2}) gives
\begin{eqnarray}\label{eqn3}
a_{2}-a_{0}\mu_{0}^{2} &=& - b_{2} \cos (\mu_{0}\tau) - b_{1} \mu_{0} \sin (\mu_{0}\tau)\\
\label{eqn4}
- \mu_{0}^{3}+a_{1}\mu_{0} &>& b_{2} \sin (\mu_{0}\tau) - b_{1} \mu_{0} \cos (\mu_{0}\tau)
\end{eqnarray}
The sufficient conditions that guarantee that stability is obtained, if (\ref{eqn3}) and (\ref{eqn4}) are satisfied simultaneously, which we utilize to get an estimate of the length of the delay. Equation (\ref{eqn3}) gives
\begin{equation}\label{eeqn6}
a_{0}\mu_{0}^{2} = a_{2}+ b_{2} \cos (\mu_{0}\tau) + b_{1} \mu_{0} \sin (\mu_{0}\tau)
\end{equation}
Taking the maximum value of $a_{2}+ b_{2} \cos (\mu_{0}\tau) + b_{1} \mu_{0} \sin (\mu_{0}\tau)$ subject to $|\sin(\mu_0\tau)|\le1, |\cos(\mu_0\tau)|\le1$, we get
\begin{eqnarray}
\nonumber a_{0} \mu_0^{2} \leq |a_2| + |b_2|+ |b_{1}| \mu_{0}\\
\label{eqn5}
\Rightarrow~~  a_{0} \mu_{0}^{2}-|b_{1}| \mu_{0} - (|a_2| + |b_2|) = 0
\end{eqnarray}
So, if
\begin{eqnarray*}
\mu_{+}=\frac{|b_{1}|+\sqrt{b_{1}^{2}+4 a_{0}(|a_{2}|+|b_{2}|)}}{2 a_{0}},
\end{eqnarray*}
then from (\ref{eqn5}), we have $\mu_{0}\leq \mu_{+}.$

From inequality (\ref{eqn4}), we get
\begin{equation}\label{eqn6}
\mu_{0}^{2}\leq a_{1}-\frac{b_{2}}{\mu_{0}} \sin(\mu_0\tau)+b_{1} \cos(\mu_0\tau),
\end{equation}
which continue to hold for sufficiently small $\tau > 0$. Substituting (\ref{eeqn6}) in (\ref{eqn6}) and rearranging we get,
\begin{equation}\label{eqn7}
~~~~~~(b_{2}-a_{0} b_{1}) (\cos(\mu_0\tau)-1)+ \left(b_{1} \mu_{0} +\frac{a_0 b_{2}}{\mu_{0}}\right) \sin(\mu_{0}\tau)< a_{0} a_{1} - a_2- b_{2}-a_{0} b_{1}
\end{equation}
Using the bounds
\nonumber \begin{equation}\label{eqnn7}
(b_{2}-a_{0} b_{1}) (\cos(\mu_0\tau)-1) = 2(b_{2}-a_{0} b_{1})\sin^{2}\left( \frac{\mu_{0} \tau}{2}\right) \leq \frac{1}{2}|b_{2}-a_{0} b_{1}| \mu_{+}^{2}\tau^{2}
\end{equation}
and
\begin{equation}\label{eqnn8}
\nonumber \left(b_{1} \mu_{0} +\frac{a_0 b_{2}}{\mu_{0}} \right) \sin(\mu_{0}\tau) \leq (b_1 \mu_{+}^{2} +a_0 b_2) \tau
\end{equation}
we obtain from (\ref{eqn7}), $N_{1}\tau^{2}+N_{2} \tau < N_{3},$ where,
\begin{eqnarray}
\nonumber N_1 = \frac{1}{2}|b_{2}- a_{0} b_{1}| \mu_{+}^{2},~~~ N_2 = (b_{1} \mu_{+}^{2}+ a_0 b_2),~~~ N_3 = a_{0} a_{1} - a_{2}-b_{2}+a_{0} b_{1}.
\end{eqnarray}
Thus, Nyquist criterion holds true for $0 \leq \tau \leq \tau_{+}$, where $\tau_{+}=\frac{1}{2 N_{1} }\left(-N_{2}+\sqrt{N_{2}^{2}+4 N_{1}N_{3}}\right),$ which gives the maximum length of delay for maintaining the stability of the limit cycle.

\subsection{Direction and stability of the Hopf bifurcation}
In the subsection 3.3, we have obtained the condition for Hopf bifurcation at the critical value $\tau_j$ and succeeded to obtain explicitly the expression for $\tau_j$ by using the normal form and manifold theory and following the ideas of Hassard et al. \cite{Hassard81}. In the coming analysis, our initial assumption is that the system (\ref{eq2}) exhibits Hopf bifurcation, which in turn imples that $\pm \textrm{i}\omega$ are the corresponding purely imaginary roots of the characteristic equation at the infected equilibrium point $E_{2}(T^{*}, I^{*}, V_{I}^{*}, V_{NI}^{*})$.

Using the transformation $u_{1}(t) = T(t)-T^{*},~ u_{2}(t) = I(t)-I^{*},~ u_{3}(t) = V_{I}(t)-V_{I}^{*},~u_{4}(t) = V_{NI}(t)-V_{NI}^{*},~ x_{i}(t)=u_{i}(\tau t), ~(i = 1,2,3,4),~ \tau = \tau_{j}+\mu,$ the system (\ref{eq2}) is converted to a functional differential equation (FDE) in $C = C([-1,~ 0],R^{4})$
\begin{equation}\label{eqn-s1-hep}
   \frac{dx}{dt}=L_{\mu}(x_{t})+f(\mu, x_{t}),
\end{equation}
where $(x_{1}(t),x_{2}(t),x_{3}(t),x_{4}(t))^{T}\in R^{4}$ and the mapping $L_{\mu}:C\rightarrow R, f:R\times C \rightarrow R$, are respectively given by
\begin{eqnarray*}\label{eqn-s2-hep}
  L_{\mu}(\phi) &=& (\tau_{j}+\mu)\left(
  \begin{array}{cccc}
    r-d_{1}-\frac{2 r T^{*}}{T_{max}}-\frac{r I^{*}}{T_{max}}-(1-c \eta_{1})\alpha V_{I}^{*}&-\frac{r T^{*}}{T_{max}}&-(1-c \eta_{1})\alpha T^{*}&0\\
    0 & -d_{2} &0&0 \\
    0 & \left(1-\frac{\eta_{r}+\eta_{1}}{2}\right)\beta & -d_{3} & 0 \\
    0 & \left(\frac{\eta_{r}+\eta_{1}}{2}\right)\beta & 0 & -d_{3}
  \end{array}
\right)
\left(
\begin{array}{c}
\phi_{1}(0)\\
\phi_{2}(0)\\
\phi_{3}(0)\\
\phi_{4}(0)
\end{array}
\right)\\
&+&(\tau_{j}+\mu)\left(
  \begin{array}{cccc}
    0&0&0&0\\
    (1-c \eta_{1})\alpha V_{I}^{*}& 0&(1-c \eta_{1})\alpha T^{*}&0 \\
    0 & 0 &0 & 0 \\
    0 & 0 & 0 & 0
  \end{array}
\right)
\left(
    \begin{array}{c}
\phi_{1}(-1)\\
\phi_{2}(-1)\\
\phi_{3}(-1)\\
\phi_{4}(-1)
\end{array}
\right)
\end{eqnarray*}
and
\begin{eqnarray*}\label{eqn-s3-hep}
~~~~~f(\mu, \phi)=(\tau_{j}+\mu)\left(
  \begin{array}{c}
  -\frac{r}{T_{max}}\phi_{1}^{2}(0)-\frac{r}{T_{max}}\phi_{1}(0)\phi_{2}(0)-(1-c \eta_{1})\alpha \phi_{1}(0) \phi_{3}(0)\\
  (1-c \eta_{1})\alpha \phi_{1}(-1) \phi_{3}(-1)\\
  0\\
  0
  \end{array}
\right)
\end{eqnarray*}
By virtue of Riesz representation theorem, we can find a function $\eta(\theta, \mu)$ of bounded variation for $\theta \in [-1,~ 0],$ such that
\begin{equation}\label{eqn-s4-hep}
L_{\mu}(\phi)=\int_{-1}^{0}  \phi(\theta) d \eta (\theta, 0)
\end{equation}
for $\phi  \in  C.$ We choose
\begin{eqnarray}
\label{eqn-s5-hep}
\nonumber \eta(\theta,\mu)&=&(\tau_{j}+\mu)
\left(
  \begin{array}{cccc}
    F&-\frac{r T^{*}}{T_{max}}&-(1-c \eta_{1})\alpha T^{*}&0\\
    0 & -d_{2} &0&0 \\
    0 & \left(1-\frac{\eta_{r}+\eta_{1}}{2}\right)\beta & -d_{3} & 0 \\
    0 & \left(\frac{\eta_{r}+\eta_{1}}{2}\right)\beta & 0 & -d_{3}
  \end{array}
\right)\delta (\theta)\\
&-&(\tau_{j}+\mu)\left(
  \begin{array}{cccc}
    0&0&0&0\\
    (1-c \eta_{1})\alpha V_{I}^{*}& 0 &(1-c \eta_{1})\alpha T^{*}&0 \\
    0 & 0 &0 & 0 \\
    0 & 0 & 0 & 0
  \end{array}
\right) \delta(\theta + 1),
\end{eqnarray}
where $\delta$ is the Dirac delta function. Let

$$
A(\mu)\phi=
\left \{\begin{array}{ll}
    {\frac{d \phi (\theta)}{d \theta},~~~~~~~~~~~~~~~~~\theta \in [-1, 0),}\\
            {\int_{-1}^{0} d \eta (\mu, s_{1}) \phi(s_{1}), ~~~\theta=0},
            \end{array}\right.
$$
and
$$
R(\mu)\phi=
\left \{\begin{array}{ll}
    {0, ~~~~~~~~~~~~\theta \in [-1, 0),}\\
    {f(\mu, \phi),~~~~ \theta = 0,}
            \end{array}\right.
$$
for $ \phi \in C^{1}([-1,~ 0], R^{4}),$ then system (\ref{eqn-s1-hep}) is equivalent to
\begin{equation}\label{eqn-s6-hep}
\dot{x_{t}}= A(\mu)x_{t}+R(\mu) x_{t},
\end{equation}
where $x_{t}(\theta)=x(t+\theta)$ for $\theta \in [-1, 0].$ Again let
$$
A^{*}\psi(s_{1})=
\left \{\begin{array}{ll}
    -{\frac{d \psi(s_{1})}{d s_{1}}, ~~~~~~~~~~~~~~~~~~~~~s_{1} \in (-1, 0],}\\
            {\int_{-1}^{0} d \eta^{T}(t, 0) \psi(-t),~~~~~~s_{1} = 0},
            \end{array}\right.
$$
for $\psi \in C^{1}([-1, 0],(R^{4})^{*})$. We define a bilinear inner product
\begin{equation}\label{eqn-s7-hep}
~~~~\langle\psi(s_{1}), \phi(\theta)\rangle=\overline{\psi}(0) \phi (0)- \int_{-1}^{0}\int_{\xi=0}^{\theta} \overline{\psi}(\xi-\theta)d \eta(\theta)\phi(\xi)d \xi,~~~[\eta(\theta)=\eta(\theta, 0)]
\end{equation}
$A(0)$ and $A^{*}$ are called as adjoint operators. It can be easily shown that $ \pm \textrm{i} \omega_{0}\tau_{j}$ are eigenvalues of $A(0)$. Therefore, $ \pm \textrm{i} \omega_{0}\tau_{j}$ are also eigenvalues of $A^{*}$.

Let $q(\theta)=(1, a, b, c_{1})^{T}e^{\textrm{i} \omega_{0}\tau_{j}\theta}$ be the eigenvector of $A(0)$ corresponding to $\textrm{i} \omega_{0}\tau_{j},$ then $A(0)q(\theta) = \textrm{i} \omega_{0}\tau_{j} q(\theta).$ Using (\ref{eqn-s2-hep}), (\ref{eqn-s4-hep}), (\ref{eqn-s5-hep}) and the definition of $A(0)$, we obtain,
$$\tau_{j}\left(
  \begin{array}{cccc}
    \textrm{i} \omega_{0}-F& \frac{r T^{*}}{T_{max}}&(1-c \eta_{1})\alpha T^{*}&0\\
    -(1-c \eta_{1})\alpha V_{I}^{*}e^{- \textrm{i} \omega_{0}\tau_{j}} & \textrm{i} \omega_{0}+d_{2} & -(1-c \eta_{1})\alpha T^{*}e^{- \textrm{i} \omega_{0}\tau_{j}}&0 \\
    0 & -\left(1-\frac{\eta_{r}+\eta_{1}}{2}\right)\beta &  \textrm{i} \omega_{0}+d_{3} & 0\\
    0 & -\left(\frac{\eta_{r}+\eta_{1}}{2}\right)\beta & 0 &  \textrm{i} \omega_{0}+d_{3}
  \end{array}
\right)q(0)=
\left(
    \begin{array}{c}
0\\
0\\
0\\
0
\end{array}
\right)$$
We can easily obtain $q(0)=(1, a, b, c_{1})^{T},$ where
$$a=\frac{( \textrm{i} \omega_{0}+d_{3})(1-c \eta_{1})\alpha V_{I}^{*}e^{- \textrm{i} \omega_{0}\tau_{j}}}{( \textrm{i} \omega_{0}+d_{2})( \textrm{i} \omega_{0}+d_{3})-(1-c \eta_{1})\left(1-\frac{\eta_{r}+\eta_{1}}{2}\right)\alpha \beta T^{*}e^{- \textrm{i} \omega_{0}\tau_{j}}},~~~ b=\frac{\left(1-\frac{\eta_{r}+\eta_{1}}{2}\right)\beta a}{( \textrm{i} \omega_{0}+d_{3})},~~~
c_{1}=\frac{\left(\frac{\eta_{r}+\eta_{1}}{2}\right)\beta a}{( \textrm{i} \omega_{0}+d_{3})}.$$
Again, let $q^{*}(s_{1})=D(1, a^{*}, b^{*}, c^{*}_{1})e^{ \textrm{i} \omega_{0}\tau_{j} s_{1}}$ be the eigenvector of $A^{*}$ corresponding to $- \textrm{i} \omega_{0}\tau_{j}.$ In the similar manner we obtain,
\begin{eqnarray}
\nonumber a^{*}=\frac{-( \textrm{i} \omega_{0}+F)}{(1-c \eta_{1})\alpha V_{I}^{*}e^{-\textrm{i} \omega_{0}\tau_{j}}},~~~~
b^{*}=\frac{\frac{r T^{*}}{T_{max}}+(-\textrm{i} \omega_{0}+d_{2})a^{*}}{\left(1-\frac{\eta_{r}+\eta_{1}}{2}\right)\beta},~~~~ c^{*}_{1}=0.
\end{eqnarray}
To assure $\langle q^{*},q(\theta) \rangle = 1,$ we use (\ref{eqn-s7-hep}) to calculate $D$ as follows:
\begin{eqnarray}
 \nonumber  \langle q^{*},q(\theta) \rangle  &=& \overline{D}(1, \overline{a^{*}}, \overline{b^{*}}, \overline{c^{*}_{1}})(1, a, b, c_{1})^{T} - \int_{-1}^{0}\int_{\xi=0}^{\theta} \overline{D}(1, \overline{a^{*}}, \overline{b^{*}}, \overline{c^{*}_{1}})e^{-i \omega_{0}\tau_{j}(\xi-\theta)}d\eta(\theta)(1, a, b, c_{1})^{T}e^{i \omega_{0}\tau_{j}\xi} d\xi\\
   \nonumber &=& \overline{D}\left\{1+a \overline{a^{*}}+b \overline{b^{*}}-\int_{-1}^{0}(1, \overline{a^{*}}, \overline{b^{*}},0)\theta e^{i\omega_{0}\tau_{j}\theta} d\eta(\theta) (1, a, b, c_{1})^{T}\right\} \\
  \nonumber  &=& \overline{D}\left\{1+a \overline{a^{*}}+b \overline{b^{*}}+\tau_{j} e^{-i \omega_{0}\tau_{j}}(1-c \eta_{1})\alpha \overline{a^{*}}(V_{I}^{*} + b T^{*}) \right\}
   \end{eqnarray}
Choosing
\begin{eqnarray}
\nonumber D=\frac{1}{1+\overline{a}a^{*}+\overline{b}b^{*}+ \tau_{j}e^{i \omega_{0}\tau_{j}}(1-c \eta_{1})\alpha a^{*}(V_{I}^{*} + b T^{*})},
\end{eqnarray}
we achieve the property $\langle q^{*},q(\theta) \rangle = 1$.

To compute the coordinates describing the center manifold  $\textbf{C}_{0}$ at $\mu =0$, we apply the idea of Hassard et al. \cite{Hassard81}. Denoting $x_{t}$ to be the solution of (\ref{eqn-s6-hep}) at $\mu =0$, we define
\begin{equation}\label{eqn-s8-hep}
z(t)=\langle q^{*},x_{t}\rangle, ~~~  W(t, \theta)=x_{t}(\theta)-2 \textrm{Re} \{z(t) q(\theta)\}
\end{equation}
On center manifold  $\textbf{C}_{0}$, we have
$$W(t, \theta)=W(z(t), \overline{z}(t),\theta),$$
where
\begin{equation}\label{eqn-s9-hep}
   W(z, \overline{z},\theta)=W_{20}(\theta)\frac{z^{2}}{2}+W_{11}(\theta)z \overline{z}+W_{02}(\theta)\frac{\overline{z}^{2}}{2}+W_{30}(\theta)\frac{z^{3}}{6}+...,
\end{equation}
z and $\overline{z}$  are local coordinates for center manifold  $\textbf{C}_{0}$ in the direction of $q^{*}$ and $\overline{q^{*}}.$ The expression $W(z, \overline{z},\theta)$ will be real provided $x_t$ is real and we are interested in real solutions only. For solution $x_{t} \in \textbf{C}_{0}$ of (\ref{eqn-s6-hep}), since $\mu = 0,$ we have
\begin{eqnarray}
   \nonumber \dot{z}(t)=\textrm{i} \omega_{0}\tau_{j} z+\bar{q^{*}}(0)f(0,W(z, \overline{z},0)+2 \textrm{Re}\{z q(\theta)\})\doteq \textrm{i} \omega_{0}\tau_{j}z+\bar{q^{*}}(0)f_{0}(z,\overline{z}).
\end{eqnarray}
which is rewritten as
\begin{eqnarray}
    \nonumber \dot{z}(t)=\textrm{i} \omega_{0}\tau_{j} z(t)+g(z,\overline{z}),
\end{eqnarray}
where
\begin{equation}\label{eqn-s10-hep}
   g(z, \overline{z})=\bar{q^{*}}(0)f_{0}(z,\overline{z})=g_{20}\frac{z^{2}}{2}+g_{11}z \overline{z} +g_{02}\frac{\overline{z}^{2}}{2} + g_{21} \frac{z^{2}\overline{z}}{2}+...
\end{equation}
Then (\ref{eqn-s8-hep}) and (\ref{eqn-s9-hep}) gives
  \begin{eqnarray*}\label{eqn-s11-hep}
  \nonumber x_{t}(\theta) &=& W(t, \theta)+2 \textrm{Re}\{z(t)q(t)\}\\
  ~~~~~~&=& W_{20}(\theta)\frac{z^{2}}{2}+W_{11}(\theta)z \overline{z}+ W_{02}(\theta)\frac{\overline{z}^{2}}{2}+(1, a, b, c_{1})^{T}e^{\textrm{i} \omega_{0}\tau_{j}\theta }z\\
  \nonumber &+& (1, \overline{a},\overline{b}, \overline{c_{1}})^{T}e^{-\textrm{i} \omega_{0}\tau_{j}\theta }\overline{z}+...,
\end{eqnarray*}
From the definition of $g(z, \overline{z})$ we have

\begin{eqnarray}\label{eqn-s12-hep}
\nonumber g(z, \overline{z})&=&\bar{q^{*}}(0)f_{0}(z,\overline{z})=\bar{q^{*}}(0)f(0, x_{t})\\
                      \nonumber &=&\tau_{j}\overline{D}(1, \overline{a^{*}}, \overline{b^{*}},0)\left(\begin{array}{c}
  -\frac{r}{T_{max}}x_{1t}^{2}(0)-\frac{r}{T_{max}}x_{1t}(0)x_{2t}(0)-(1-c \eta_{1})\alpha x_{1t}(0) x_{3t}(0)\\
  (1-c \eta_{1})\alpha x_{1t}(-1) x_{3t}(-1)\\0\\0\\ \end{array} \right)\\
\nonumber &=& -\tau_{j}\overline{D}\{\frac{r}{T_{max}} \left[z+\overline{z}+W_{20}^{(1)}(0)\frac{z^{2}}{2}+W_{11}^{(1)}(0)z \overline{z}+W_{02}^{(1)}(0)\frac{\overline{z}^{2}}{2}+ o(|(z, \overline{z})|^{3})  \right]^{2}\\
   \nonumber &+& \frac{r}{T_{max}} \left[z+\overline{z}+W_{20}^{(1)}(0)\frac{z^{2}}{2}+W_{11}^{(1)}(0)z \overline{z}+W_{02}^{(1)}(0)\frac{\overline{z}^{2}}{2}+ o(|(z, \overline{z})|^{3})  \right]\\
   && ~~~~~~~~~~~\times \left[a z+ \bar{a}\overline{z}+W_{20}^{(2)}(0)\frac{z^{2}}{2}+W_{11}^{(2)}(0)z \overline{z}+W_{02}^{(2)}(0)\frac{\overline{z}^{2}}{2}+ o(|(z, \overline{z})|^{3})  \right] \\
   \nonumber  &+&\alpha_{1}\left[ b e^{-\textrm{i} \omega_{0}\tau_{j}}z+\overline{b} e^{\textrm{i} \omega_{0}\tau_{j}}\overline{z}+W_{20}^{(3)}(-1)\frac{z^{2}}{2}+W_{11}^{(3)}(-1)z \overline{z}+W_{02}^{(3)}(-1)\frac{\overline{z}^{2}}{2}+ o(|(z, \overline{z})|^{3}) \right]\\
   \nonumber &\times& \left[e^{-\textrm{i} \omega_{0}\tau_{j}}z+e^{\textrm{i} \omega_{0}\tau_{j}}\overline{z}+W_{20}^{(1)}(-1)\frac{z^{2}}{2}+W_{11}^{(1)}(-1)z \overline{z}+W_{02}^{(1)}(-1)\frac{\overline{z}^{2}}{2}+ o(|(z, \overline{z})|^{3}) \right]\}\\
   \nonumber  &+& \gamma \left[ b e^{-\textrm{i} \omega_{0}\tau_{j}}z+\overline{b} e^{\textrm{i} \omega_{0}\tau_{j}}\overline{z}+W_{20}^{(3)}(-1)\frac{z^{2}}{2}+W_{11}^{(3)}(-1)z \overline{z}+W_{02}^{(3)}(-1)\frac{\overline{z}^{2}}{2}+ o(|(z, \overline{z})|^{3}) \right]\\
   \nonumber  &\times& \left[e^{-\textrm{i} \omega_{0}\tau_{j}}z+e^{\textrm{i} \omega_{0}\tau_{j}}\overline{z}+W_{20}^{(1)}(-1)\frac{z^{2}}{2}+W_{11}^{(1)}(-1)z \overline{z}+W_{02}^{(1)}(-1)\frac{\overline{z}^{2}}{2}+ o(|(z, \overline{z})|^{3}) \right]
\end{eqnarray}
Comparing the coefficients with (\ref{eqn-s10-hep}), we get,
\begin{eqnarray*}\label{eqn-s13-hep}
~~~~~~~~\nonumber  g_{20} &=& -2 \tau_{j}\overline{D}\left[(1+a)\frac{r}{T_{max}}+(1-c\eta_{1})\alpha b e^{-2 \textrm{i} \omega_{0}\tau_{j}}-(1-c\eta_{1})\alpha b \overline{a^{*}}e^{-2 \textrm{i} \omega_{0}\tau_{j}}  \right]\\
\nonumber  g_{11} &=& -2 \tau_{j}\overline{D}\left[(1+Re(a))\frac{r}{T_{max}}+(1-c\eta_{1})\alpha Re(b)-(1-c\eta_{1})\alpha Re(b) \overline{a^{*}}\right] \\
  g_{02} &=& -2 \tau_{j}\overline{D}\left[(1+\overline{a})\frac{r}{T_{max}}+(1-c\eta_{1})\alpha \overline{b} e^{2 \textrm{i} \omega_{0}\tau_{j}}-(1-c\eta_{1})\alpha \overline{b} \overline{a^{*}}e^{2 \textrm{i} \omega_{0}\tau_{j}}  \right] \\
\nonumber  g_{21} &=& -2 \tau_{j}\overline{D}[2 \frac{r}{T_{max}}W_{11}^{(1)}(0)+ \frac{r}{T_{max}}W_{20}^{(1)}(0)+ \frac{r}{T_{max}}\{a W_{11}^{(1)}(0)+\frac{1}{2}\overline{a} W_{20}^{(1)}(0)\\
\nonumber &+&\frac{1}{2}W_{20}^{(2)}(0)+W_{11}^{(2)}(0) \}+(1-c\eta_{1})\alpha \{W_{11}^{(3)}{(-1)}e^{-\textrm{i} \omega_{0}\tau_{j}}+\frac{1}{2}W_{20}^{(3)}{(-1)}e^{\textrm{i} \omega_{0}\tau_{j}} \\
\nonumber &+& \frac{1}{2}\overline{b}W_{20}^{(1)}{(-1)}e^{\textrm{i} \omega_{0}\tau_{j}}+b W_{11}^{(1)}{(-1)}e^{-\textrm{i} \omega_{0}\tau_{j}}\}-(1-c\eta_{1})\alpha \overline{a^{*}}\{W_{11}^{(3)}{(-1)}e^{-\textrm{i} \omega_{0}\tau_{j}} \\
\nonumber &+& \frac{1}{2}W_{20}^{(3)}{(-1)}e^{\textrm{i} \omega_{0}\tau_{j}}+\frac{1}{2}\overline{b} W_{20}^{(1)}{(-1)}e^{\textrm{i} \omega_{0}\tau_{j}}+ b W_{11}^{(1)}{(-1)}e^{-\textrm{i} \omega_{0}\tau_{j}}\} ]
\end{eqnarray*}
Since $W_{20}(\theta)$ and $W_{11}(\theta)$ are in the expression of $g_{21}$, we need to compute them. The calculation is shown in the Appendix.
Furthermore, $g_{21}$ can be expressed using the parameters and the delay in the system. Therefore, we can compute
\begin{eqnarray}\label{bifur}
\nonumber c_{11}(0) &=& \frac{\textrm{i}}{2 \omega_{0}\tau_{j}}\left(g_{20}g_{11}-2|g_{11}|^{2}-\frac{|g_{02}|^{2}}{3}\right)+\frac{g_{21}}{2}, \\
\nonumber  \mu_{2}&=& -\frac{\textrm{Re}\{c_{11}(0)\}}{\textrm{Re}\{\lambda^{'}(\tau_{j})\}}, \\
\beta_{2}&=& 2 \textrm{Re}(c_{11}(0))~~ \textrm{and} \\
\nonumber T_{2}&=& -\frac{\textrm{Im} \{ c_{11}(0)\}+ \mu_{2}  \textrm{Im}\{\lambda^{'}(\tau_{j})\}}{\omega_{0}\tau_{j}}.
\end{eqnarray}
The expressions in (\ref{bifur}) gives the quantification of the bifurcating periodic solutions in the center manifold at the critical value ($\tau_0$) of the time lag. The direction of the Hopf bifurcation is determined from the sign of $\mu_2$. The Hopf bifurcation is forward or backward according as $\mu_2 > 0$ or $\mu_2 < 0$ and the bifurcating periodic solution exists for $\tau > \tau_0$ and $\tau > \tau_0$ respectively. The sign of $\beta_2$ quantifies the stability of the bifurcating periodic solution; for $\beta_2 < 0$, the solution is stable and unstable if $\beta_2 > 0$. Finally, $T_2$ determines the period of the bifurcating periodic solution. The period is increasing if $T_2 > 0$ and decreasing if $T_2 < 0$.

\section{Numerical results and its biological implications}
We now simulate the system (\ref{eq2}) to observe the effect of intracellular delay on the dynamics of hepatitis C virus with therapy. We verify our model with two sets of parameter obtained from \cite{Dahari07} (see table 1 and 2), starting with parameter set 1. We have taken the viral replication time $(\tau)$ in the range of $24$ hours to $36$ hours \cite{NarendraM04, Guatelli90}. Throughout the discussion, we have considered $\eta$ as the combined drug efficacy, given by
\begin{eqnarray*}
1-\eta = (1-c \eta_{1})\left(1-\frac{\eta_{r}+\eta_{1}}{2}\right),
\end{eqnarray*}
along with $\eta_c$, the critical drug efficacy. Our discussion will revolve around $R_0$, the basic reproduction number and its relation with the critical drug efficacy $\eta_c$. Please note, the initial viral load for numerical simulation is $10^7$ copies/ml, which is around the endemic equilibrium point and hence
\begin{eqnarray}
\nonumber R_0 > 1 ~~\Rightarrow~~\frac{\hat{T}}{T^*} > 1 ~~\Rightarrow~~\hat{T} \frac{\alpha \beta (1-\eta)}{d_2 d_3} > 1 \\
\Rightarrow~~\eta < 1-\frac{d_2 d_3}{\alpha \beta \hat{T}} = \eta_c
\end{eqnarray}
This implies when $\eta < \eta_c$, $R_0 > 1$ and the viral load reaches the endemic equilibrium point. But, when $\eta > \eta_c$, $R_0 < 1$, the endemic equilibrium point becomes negative and exchanges stability with the uninfected steady state and disease free equilibrium point is achieved.

~~~~~~~~~~~~~~~~~~~~~~~~~~~~~~~~~~~~~~~~~~~~\textbf{Table 1}
\begin{eqnarray}
 \nonumber \begin{tabular}{|c|c|c|}
\hline
Parameter&Parameter values& Units\\
\hline
$s$ & $1.0$  & $\textmd{cells}~\textmd{day}^{-1}~\textmd{ml}^{-1} \cite{Dahari07}$\\
$r$ & $2.0$  & $\textmd{day}^{-1} \cite{Dahari07} $\\
$T_{\max}$ & $3.6 \times 10^{7}$  & $\textmd{cells}~\textmd{ml}^{-1} \cite{Dahari07}$ \\
$\alpha$ & $2.25 \times 10^{-7}$  & $\textmd{ml}~\textmd{virions}^{-1}~\textmd{day}^{-1} \cite{Dahari07} $\\
$\beta$ & $2.9$  & $\textmd{virions}~\textmd{day}^{-1} \cite{Dahari07}$\\
$d_1$ & $0.01$  & $\textmd{day}^{-1} \cite{Dahari07}$ \\
$d_2$ & $1.0$  & $\textmd{day}^{-1} \cite{Dahari07}$ \\
$d_3$ & $6.0$  & $\textmd{day}^{-1} \cite{Dahari07}$\\
%$c$ & $0.08$  & dimensionless\\
%$\eta_{1}$ & $0.5$  & dimensionless\\
%$\eta_{r}$ & $0.7$  & dimensionless\\
\hline
\end{tabular}
\end{eqnarray}

%\textbf{Case 1} We have a set of parameter values given in table 1. For this parameter values with $c = 0.08$, $\eta_{1}= 0.5$, $\eta_{r} = 0.7$, $E_{2}= (2.39464 \times 10^{7}, 6.77843 \times 10^{6}, 1.3105 \times 10^{6}, 1.96574 \times 10^{6})$. We can not find any positive value of $\omega$ and $\tau$ such that transversality condition of Hopf bifurcation $\left[\frac{d(Re \lambda)}{d\tau}\right]_{\omega~=~\omega_0, \tau~=~\tau_j}>0$ is holds. It means for this set of parameter values Hopf bifurcation does not exist.
%Now, we have taken $\tau = 5.0$ for set of parameter values given in table 1 with different values of $c, \eta_{1}, \eta_{r}$, get different dynamics of HCV which have shown in fig 7.1 to fig 7.6. fig 7.7 shows that we have fitted HCV dynamics with different patients data taken from \cite{dahari2010} and different values of $\tau$

~~~~~~~~~~~~~~~~~~~~~~~~~~~~~~~~~~~~~~~~~~~~\textbf{Table 2}
\begin{eqnarray}
 \nonumber \begin{tabular}{|c|c|c|}
\hline
Parameter&Parameter values& Units\\
\hline
$s$ & $3.7 \times 10^{4}$  & $\textmd{cells}~\textmd{day}^{-1}~\textmd{ml}^{-1}$ \cite{Dahari07}\\
$r$ & $0.73$  & $\textmd{day}^{-1}$ \cite{Dahari07}\\
$T_{\max}$ & $0.6 \times 10^{7}$  & $\textmd{cells}~\textmd{ml}^{-1}$ \cite{Dahari07}\\
$\alpha$ & $1.8 \times 10^{-7}$  & $\textmd{ml}~\textmd{virions}^{-1}~\textmd{day}^{-1}$ \cite{Dahari07}\\
$\beta$ & $13.9$  & $\textmd{virions}~\textmd{day}^{-1}$ \cite{Dahari07}\\
$d_1$ & $2.4 \times 10^{-3}$  & $\textmd{day}^{-1}$ \cite{Dahari07}\\
$d_2$ & $0.06$  & $\textmd{day}^{-1}$ \cite{Dahari07}\\
$d_3$ & $13.9$  & $\textmd{day}^{-1}$ \cite{Dahari07}\\
%$c$ & $0.08$  & dimensionless\\
%$\eta_{1}$ & $0.5$  & dimensionless\\
%$\eta_{r}$ & $0.7$  & dimensionless\\
\hline
\end{tabular}
\end{eqnarray}

\begin{figure}
\centering
\includegraphics[width=4.0in]{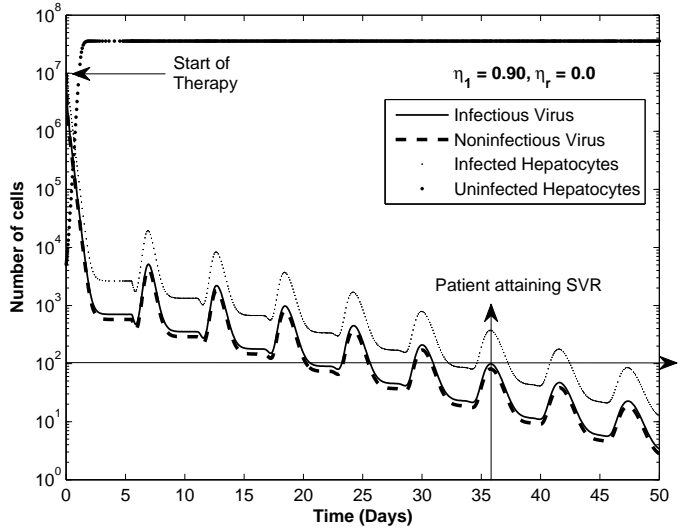}
\caption{\emph{Interferon monotherapy in presence of intracellular delay ($\tau = 22.0$ hours), where the parameter set is taken from Table 1. with initial conditions $(10^7, 10^7, 10^7, 10^7)$. The dynamics of all the state variables is shown in the figure. The hepatitis C viral load declines, which is biphasic in nature. From the start of the therapy, the first phase of virus decline occurs for (48-52) hours followed by a slow second phase till the patient reaches SVR in 36 days (shown by intersection of arrows). The infected hepatocytes and non-infectious virus also follow the same pattern of decline. The wavy nature of the decline is due to the presence of intracellular delay. The uninfected hepatocytes reaches a high steady state value ranging between $(10^7 - 10^8)$ cells, which has been observed in healthy patients. Please note that the uninfected hepatocytes though starts with the initial value of $10^7$ cells, sharply declines when the viral load is high initially and regains its normal range once SVR has been achieved.}}
\label{fig1}
\end{figure}

\textbf{Effect of interferon monotherapy:} We first consider interferon, a mixture of proteins, having antiviral and immunomodulating effects, to illustrates interferon monotherapy on viral dynamics in the presence of intracellular delay. Hepatitis C virus declines as well as the dynamics of infected and uninfected hepatocytes are shown in figure \ref{fig1}.

\begin{figure}
\centering
\includegraphics[width=4.0in]{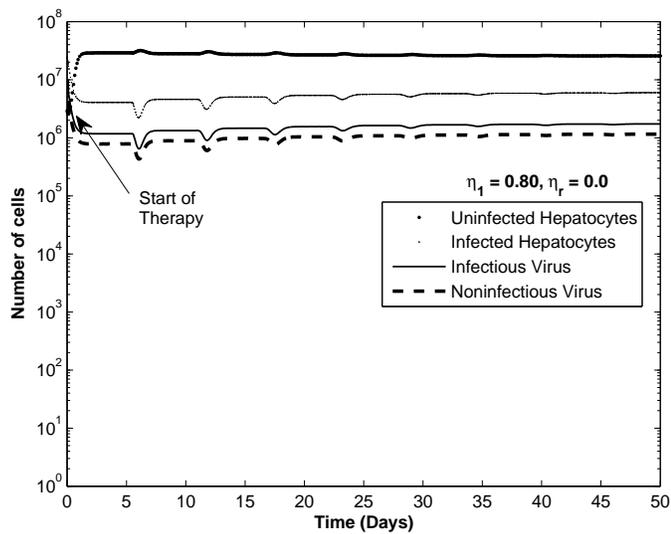}
\caption{\emph{Non-responders in the case of Interferon monotherapy in presence of intracellular delay ($\tau = 22.0$ hours), where the parameter set is taken from Table 1. with initial conditions $(10^7, 10^7, 10^7, 10^7)$. There is a short first phase decline, followed by a flat phase, where the hepatitis C virus ceases to decline and reaches the endemic steady state. This means that the drug interferon fails to reach the critical drug efficacy, so needed to eradicate hepatitis C virus and achieve SVR.}}
\label{fig2}
\end{figure}

The viral decline is biphasic in nature, mimicking the kinetics of the viral decline in patients responding to interferon monotherapy. From the figure we conclude that the first phase decline occurs for about (48-52) hours, following by a slow second phase decline starting  at 52 hours. This behavior has been observed clinically \cite{Lam97, Zeuzem98}. The wavy nature of the second phase decline is due to the effect of intracellular delay. As the viral load declines, the infected hepatocytes also decline and the number of uninfected hepaocytes reaches large steady  state value. However, the standard interferon regime sometimes fails to control HCV, which are the case of non-responders. This dynamics has been portrayed in figure \ref{fig2}, where there is a first phase decline followed by a flat phase, that is, the viral load approaches an endemic steady state. The infectious hepatocytes and noninfectious viruses also exhibit the same behavior. Please note, though the interferon efficacy is high ($\eta_1 = 0.80$), the combined efficacy $\eta = 0.64 < \eta_c (=0.745)$. Since, $\eta < \eta_c$ implies $R_0 > 1$, the viral load showing minor decline will never converge to $E_0$ and ultimately converge to a new steady state. This means that in non-responders, the efficacies due to interferon fail to reach the critical drug efficacy, and this has been confirmed clinically too \cite{Bekkering97, Bekkering98}.\\

\begin{figure}
\centering
\includegraphics[width=4.0in]{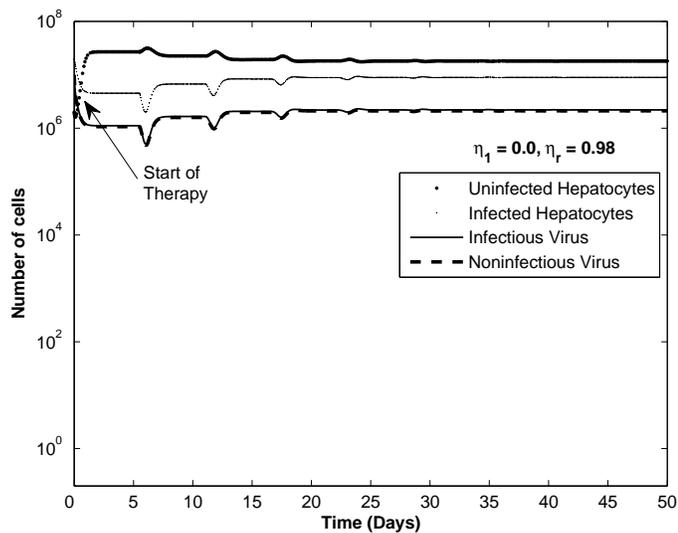}
\caption{\emph{Ribavirin monotherapy in presence of intracellular delay ($\tau = 22.0$ hours), where the parameter set is taken from Table 1. with initial conditions $(10^7, 10^7, 10^7, 10^7)$. The dynamics of all the state variables is shown in the figure, where ribavirin fails to clear the viral load and hepatitis C virus remains more or less constant during the therapy.}}
\label{fig3}
\end{figure}

\textbf{Effect of ribavirin monotherapy:} Ribavirin monotherapy has not been effective for the treatment of hepatitis C viral infection and our model also certifies that. The effect of ribavirin monotherapy on HCV is shown in figure \ref{fig3}, where ribavirin fails to clear the viral load. Even though the efficacy of ribavirin is 0.98 (figure \ref{fig3}), its combined efficacy as calculated by is 0.68, which is less than $\eta_c (=0.745)$, implying $R_0 > 1$. Several experiments have been conducted \cite{Di Bisceglie92, Di Bisceglie95, Bodenheimer97} to evaluate the effect of ribavirin for the treatment of hepatitis C virus but no virologic end of treatment response (ETR) were observed in all these experiments. In most of these cases, HCV RNA level remains constant after 3-6 months of ribavirin monotherapy compared within pretreatment values \cite{Di Bisceglie92, Di Bisceglie95}.\\

\begin{figure}
\centering
\includegraphics[width=4.0in]{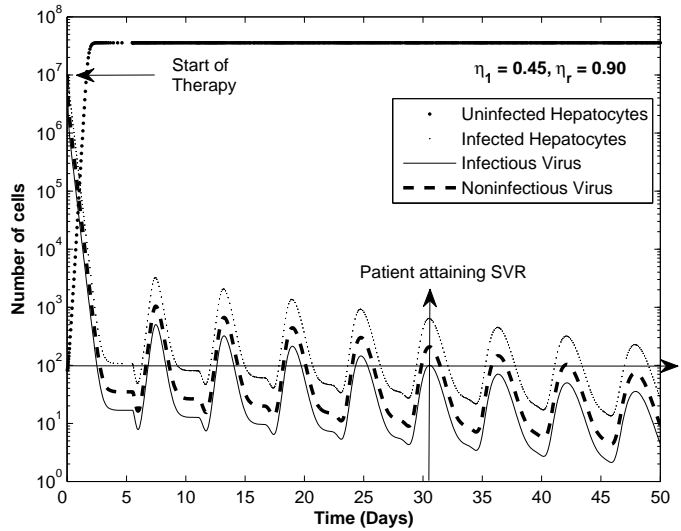}
\caption{\emph{Dynamics of hepatitis C virus during combined therapy in presence of intracellular delay ($\tau = 22.0$ hours), where the parameter set is taken from Table 1 with initial conditions $(10^7, 10^7, 10^7, 10^7)$. From the figure it is observed that when interferon efficacy is small, ribavirin enhances treatment level by turning infectious virus to non-infectious and the concentration of virus in the patient's blood falls below the detection limit of 100 copies/ml. The pattern of decline is wavy due to the delay effect and the patient attains SVR in 31 days (shown by intersection of arrows), which is a little faster than interferon monotherapy.}}
\label{fig4}
\end{figure}

\begin{figure}
\centering
\includegraphics[width=4.0in]{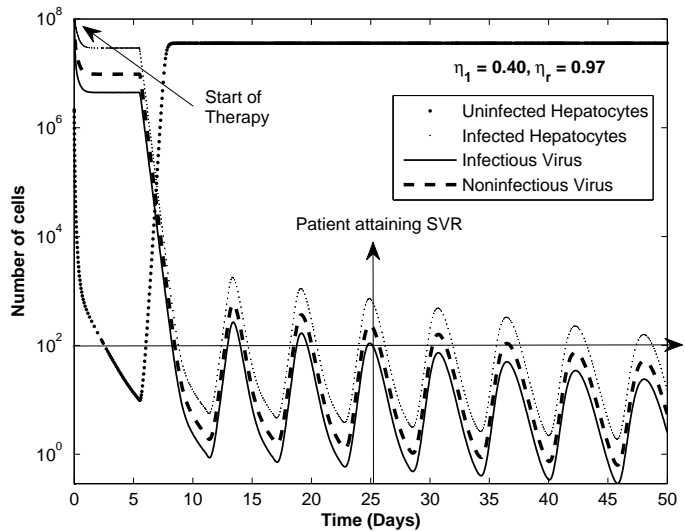}
\caption{\emph{Triphasic decline in hepatitis C viral dynamics in presence of intracellular delay ($\tau = 22.0$ hours), where the parameter set is taken from Table 1 with initial conditions $(10^7, 10^7, 10^7, 10^7)$. The first phase decline occurs for (24-36) hours followed by a shoulder (flat) phase, which continues for (1.5-5) days. Thus, there is a third phase decline, which happens in oscillation due to the presence of delay and SVR is achieved in 25 days.}}
\label{fig5}
\end{figure}

\textbf{Effect of combined therapy:} To achieve the goal of minimizing the viral load of HCV, current standard treatment consists of the combination of pegylated interferon and ribavirin \cite{Layden03}. Thus, we conclude that when interferon efficacy is small, ribavirin enhances treatment level by turning progeny virions non-infections, resulting in a second phase decline and the patient attains SVR, that is, the concentration of the virus in the patient's blood falls below the detection limit (100 copies/ml) during the combined therapy (figure \ref{fig4}). Experimental observation also confirm this fact \cite{Herrmann03, Layden03, Pawlotsky04}. The combined therapy of pegylated interferon and ribavirin also results in triphasic decline of viral load in some patients \cite{Herrmann03}. Figure \ref{fig5} shows triphasic decline where the first phase shows a sharp drop in the viral load followed by a shoulder second phase (flat phase) after which the viral load decline continues in oscillations (third phase decline) and the patient attains SVR in approximately 25 days.\\

\begin{figure}
\centering
\includegraphics[width=6.0in]{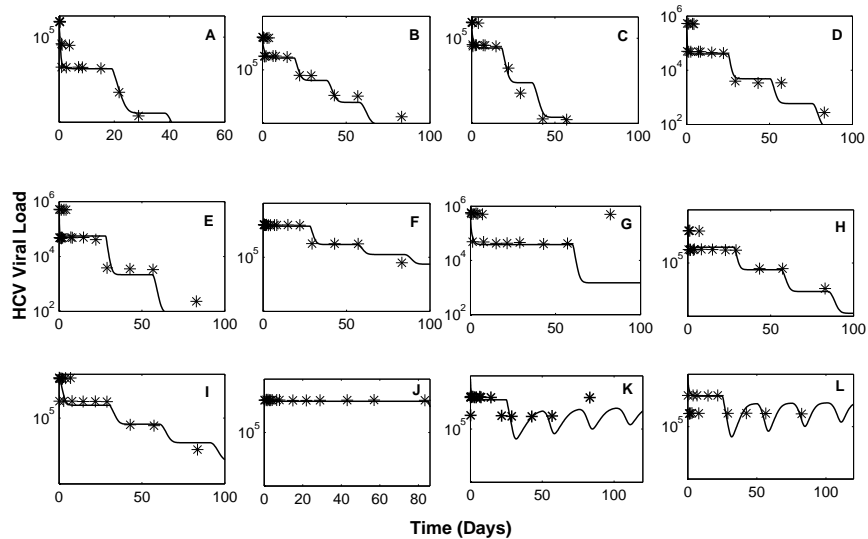}
\caption{\emph{Model validation with data obtained from 12 patients who were undergoing treatment with PEG-IFN-$\alpha$-2a serum concentration at the University of Sao Paulo hospital das clinicas. The hepatitis C virus behavioral data of the patients which are shown in $*$ are obtained from \cite{dahari2010}, page 464, figure2. In this figure, A-I are the responses of the model with interferon monotherapy (continuous line) which matches with the data ($*$). J is non-responder, K and L are non-SVR. All the parameter values are taken from Table 1, with $tau = 22$ hours.}}
\label{fig6}
\end{figure}

Figure \ref{fig6} shows that the dynamics of our model are in agreement with the behavioral pattern of HCV in 12 patients obtained from the study, conducted at the University of Sao Paulo Hospital das clinicas and was approved by the ethics in research committee \cite{dahari2010}. The patients were treated with pegylated interferon and ribavirin for 48 weeks and their detailed analysis were recorded (denoted by $*$ in the figure). Figure A, B, C, D, E, F, G, H, I show the behavior of HCV RNA levels during the first 12-14 weeks of PEG-IFN-$\alpha$-2a monotherapy. The graph obtained from the model (continuous line) matches with the patient data ($*$) \cite{dahari2010}. Figure J is the case of non-responders and figures K and L are the case of non-SVR.\\

\begin{figure}
\centering
\includegraphics[width=7in]{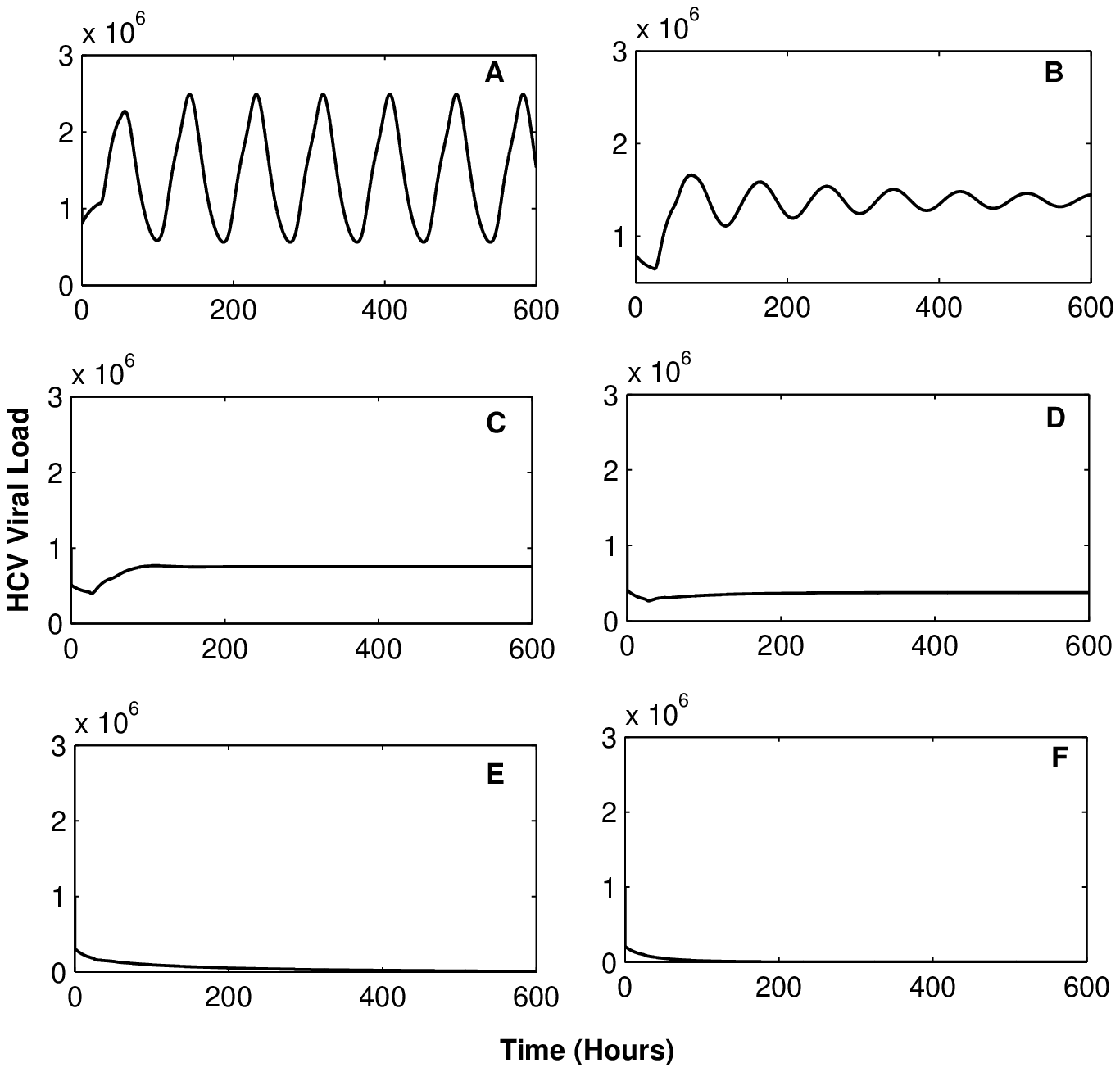}
\caption{\emph{The figure shows the dynamics of hepatitis C virus with intracellular delay $\tau$ = 24.1 hours (= $\tau_0$), $\eta_1$ = 0.5, $\eta_r$ = 0.7, c = 0.81, where the system exhibits Hopf bifurcation. As $\eta_r$ is increased to 0.75, 0.8, 0.85, 0.9 and 0.95, decline in viral load has been observed (B, C, D, E, F). All the parameter values are taken from Table 2.}}
\label{fig7}
\end{figure}

We next consider the second set of parameters given in the table 2, for verification of the model (\ref{eq2}). With these set of parameters, theorem 3.1 holds. It is observed that without delay there exist a unique interior equilibrium point  $E_{2}$, which is locally asymptotically stable. Using the same set of parameters, we obtain, we $A_{1} > 0$  and $A_{3} < 0$, which indicates the existence of a positive root. We also calculate $w_{0} = 0.07321$ by solving (\ref{transcendentaleq3}). Following theorem 1 in \cite{Cooke86}, we can say that stability switch may occur as the value of $\tau$ increases. In this case stability switch occurs at $\tau_{0} = 24.1$ hours (calculated from (\ref{proposition1eq4}) by putting $j = 0$), where the dynamics of the system changes from a stable steady state to stable oscillatory state. The calculated value of  $\tau_{0}$ also certifies our analytical estimation of the length of delay to preserve stability. Therefore, by Butler's lemma, $E_{2}$ remains stable for $\tau < \tau_{0}$ (= 24.1 hours). For $\tau = 24.1$ hours Hopf  bifurcation occurs resulting in stable periodic solution and at the bifurcated point, a stable limit cycle is formed (figure \ref{fig7}). We observe a stability switch as $\tau$ crosses the threshold $(\tau = 24.1)$ hours and shows aperiodic solutions. Please note that we have not provided figures for $\tau < \tau_{0}$ and $\tau > \tau_{0}$ as these values of $\tau_0$ are unrealistic as far as the patient's data are concerned.

%
%\begin{figure}\label{fig7}
%\centering
%\includegraphics[width=7in]{delay_april10_07.eps}
%\caption{\emph{This shows that uninfected hepatocytes, infected hepatocytes, infectious virus and noninfcetious virus converge to their equilibrium point with $\tau= 3.0$ respectively. The parameter set is taken from Table 2.}}
%\end{figure}
%%
%
%\begin{figure}\label{fig9}
%\centering
%\includegraphics[width=7in]{delay_april10_09.eps}
%\caption{\emph{This shows the oscillations of uninfected hepatocytes, infected hepatocytes, infectious virus and noninfcetious virus with $\tau= 80.0$. The parameter set is taken from Table 2.}}
%\end{figure}

\section{Summary and Discussions}

Several studies have used mathematical models to provide insights into the mechanism and dynamics of the progression of hepatitis C viral infection. It is of both mathematical and biological interest to determine the cause of sustained oscillation, which can be the result of intracellular delay. The main goal of this paper is to investigate how the delay affects the overall disease progression and mathematically, how the dynamics of the system is affected by the delay.\\

We have studied a modified HCV model with four state variables, namely, target cells, infected hepatocytes, infected and non-infected virions, along with an intracellular delay and combined drug therapy (interferon and ribavirin). The basic reproduction number of the hepatitis C virus has been used largely in understanding the dynamics of the viral infection. The influence of time delay on the stability of the equilibrium states has been discussed. We have proved that the local stability of the disease free steady state is independent of the delay. However, for the endemic equilibrium state, an increase in delay can destabilize the system, leading to hopf bifurcation and periodic solutions. The estimation of the length of delay to preserve stability has been obtained. Also, the direction and stability of the Hopf bifurcation have been studied. \\

The numerical study of the system focused on two aspects, namely, the dynamics of the system on two different parameter sets and validation of the model with the virus behavioral pattern of 12 patients, under the influence of interferon monotherapy and combined therapy with ribavirin. Our results mimic the behavior of HCV in patients observed clinically with the first set of data given in the table 1. The second set of data actually calculate a discrete time delay $\tau_0$, for which we get a Hopf bifurcation resulting in stable periodic solution. It is numerically observed that with proper choice of drug efficacies, hepatitis C virus can be controlled. Please note that a Hopf bifurcation was not found in a realistic parameter space with data set 1 but with data set 2, our analytical finding on the existence of a hopf bifurcation with biologically realistic parameters remain true. In the future, we would like to extend this work for different types of delays such as distributed delays or multiple discrete delays. Also, there are many components in this model that may be regarded as stochastic rather than deterministic and these variations may significantly alter the dynamics of the system. We, therefore, suggest to convert the system given by (\ref{eq2}) into a stochastic delay differential equations and study its dynamics, which we propose as our future work.

\newpage

\section*{Appendix: Calculation of $W_{20}(\theta)$ and $W_{11}(\theta)$}

From (\ref{eqn-s6-hep}) and (\ref{eqn-s8-hep}), we have
\begin{eqnarray}\label{eqn-s14-hep}
 \dot{W}=\dot{x}_{t}-\dot{z} q-\dot{\bar{z}}\bar{q}
&=&\left \{\begin{array}{ll}
    {A W- 2 \textrm{Re}\{\bar{q^{*}}(0)f_{0} q (\theta)\},~~~~~~ if \theta \in  [-1, 0),}\\
            {A W-2 \textrm{Re}\{\bar{q^{*}}(0)f_{0} q (0)\}+f_{0},~~if \theta=0},
            \end{array}\right. \\
            \nonumber &=& A W+H(z, \bar{z}, \theta)~~~~(say),
 \end{eqnarray}
 where
\begin{equation}\label{eqn-s15-hep}
  H(z, \bar{z}, \theta)=H_{20}(\theta)\frac{z^{2}}{2}+H_{11}(\theta)z \overline{z} +H_{02}(\theta)\frac{\overline{z}^{2}}{2} +...
\end{equation}
Substituting the corresponding series into (\ref{eqn-s14-hep}) and comparing the coefficients, we obtain
\begin{equation}\label{eqn-s16-hep}
  (A-2 \textrm{i} \tau_{0}\omega_{0})W_{20}(\theta) =-H_{20}(\theta), A W_{11}(\theta)=-H_{11}(\theta),...
\end{equation}
From (\ref{eqn-s14-hep}), we know that for $\theta \in [-1, 0),$
\begin{equation}\label{eqn-s17-hep}
  H(z, \bar{z}, \theta) =-\bar{q^{*}}(0)f_{0} q (\theta)- q^{*}(0) \bar{f_{0}} \bar{q}(\theta)=-g(z, \bar{z})q(\theta)-\bar{g}(z, \bar{z})\bar{q}(\theta).
\end{equation}
 Comparing the corresponding coefficients with that of  (\ref{eqn-s15-hep}), we get,
 \begin{equation}\label{eqn-s18-hep}
  H_{20}(\theta) =-g_{20}q(\theta)-\bar{g}_{02}\bar{q}(\theta),
  H_{11}(\theta) =-g_{11}q(\theta)-\bar{g}_{11}\bar{q}(\theta).
\end{equation}
 From (\ref{eqn-s16-hep}),(\ref{eqn-s18-hep}) and the definition of $A$, it follows that
\begin{eqnarray}\label{eqn-s19-hep}
\nonumber  \dot{W_{20}}(\theta)=2 \textrm{i}\omega_{0}\tau_{k}W_{20}(\theta)+g_{20}(\theta)q(\theta)+\bar{g}_{02}\bar{q}(\theta).
\end{eqnarray}
 Note that  $q(\theta)=(1, a, b, c_{1})^{T}e^{\textrm{i} \omega_{0}\tau_{j}\theta}$, hence
  \begin{equation}\label{eqn-s20-hep}
  W_{20}(\theta)=\frac{\textrm{i}g_{20}}{\omega_{0}\tau_{j}}q(0)e^{\textrm{i} \omega_{0}\tau_{j}\theta}+\frac{\textrm{i}\bar{g}_{02}}{3\omega_{0}\tau_{j}}\bar{q}(0)e^{-\textrm{i} \omega_{0}\tau_{j}\theta}+E_{1}e^{2\textrm{i} \omega_{0}\tau_{j}\theta},
\end{equation}
Here $E_{1}=(E_{1}^{(1)}, E_{1}^{(2)}, E_{1}^{(3)}, E_{1}^{(4)})\in R^{4}$ is a constant vector. Similarly, from (\ref{eqn-s16-hep}) and (\ref{eqn-s19-hep}), we get
\begin{equation}\label{eqn-s21-hep}
  W_{11}(\theta)=-\frac{\textrm{i}g_{11}}{\omega_{0}\tau_{j}}q(0)e^{\textrm{i} \omega_{0}\tau_{j}\theta}+\frac{\textrm{i}\bar{g}_{11}}{3\omega_{0}\tau_{j}}\bar{q}(0)e^{-\textrm{i} \omega_{0}\tau_{j}\theta}+E_{2},
\end{equation}
 where $E_{2}=(E_{2}^{(1)}, E_{2}^{(2)}, E_{2}^{(3)}, E_{2}^{(4)})\in R^{4}$ is also a constant vector. Now, we have to find an appropriate constant vector $E_1$ and $E_2$ which satisfy the above conditions. From the definition of $A$ and (\ref{eqn-s16-hep}), we obtain
\begin{equation}\label{eqn-s22-hep}
\int_{-1}^{0} d \eta (\theta) W_{20}(\theta)= 2 \textrm{i}\omega_{0} \tau_{j}W_{20}(\theta)-H_{20}(\theta),
\end{equation}
and
 \begin{equation}\label{eqn-s23-hep}
\int_{-1}^{0} d \eta (\theta) W_{11}(\theta)= -H_{11}(0),
\end{equation}
where $\eta(\theta)=\eta(0, \theta).$ Using (\ref{eqn-s14-hep}), we have
\begin{eqnarray}\label{eqn-s24-hep}
~~~~~~~~H_{20}(0)=-g_{20}q(0)-\bar{g}_{20}\bar{q}(0)+2 \tau_{j}\left(
\begin{array}{c}
-(1+a)\frac{r}{T_{max}}-(1-c\eta_{1})\alpha b e^{-2 \textrm{i} \omega_{0}\tau_{j}} \\
(1-c\eta_{1})\alpha b e^{-2 \textrm{i} \omega_{0}\tau_{j}}  \\
0 \\
0 \\
\end{array}
\right)
\end{eqnarray}
and
\begin{eqnarray}\label{eqn-s25-hep}
~~~~~~~~~H_{11}(0)=-g_{11}q(0)-\bar{g}_{11}\bar{q}(0)+2 \tau_{j}\left(
\begin{array}{c}
-(1+Re(a))\frac{r}{T_{max}}-(1-c\eta_{1})\alpha Re(b) \\
(1-c\eta_{1})\alpha Re(b)  \\
0 \\
0 \\
\end{array}
\right)
\end{eqnarray}\\
Substituting (\ref{eqn-s20-hep}) and (\ref{eqn-s24-hep}) into (\ref{eqn-s22-hep}), we obtain\\\\
$\left(2 \textrm{i}\omega_{0} \tau_{j} I-\int_{-1}^{0}e^{2 \textrm{i} \omega_{0}\tau_{j}\theta} d \eta(\theta)\right)E_{1}=2 \tau_{j}\left(
\begin{array}{c}
-(1+a)\frac{r}{T_{max}}-(1-c\eta_{1})\alpha b e^{-2 \textrm{i} \omega_{0}\tau_{j}} \\
(1-c\eta_{1})\alpha b e^{-2 \textrm{i} \omega_{0}\tau_{j}}  \\
0 \\
0 \\
\end{array}
\right)
$\\
 which leads to\\\\
 $\left(
  \begin{array}{cccc}
    i \omega_{0}-F& \frac{r T^{*}}{T_{max}}&(1-c \eta_{1})\alpha T^{*}&0\\
    -(1-c \eta_{1})\alpha V_{I}^{*}e^{-i \omega_{0}\tau_{j}} &i \omega_{0}+d_{2} & -(1-c \eta_{1})\alpha T^{*}e^{-i \omega_{0}\tau_{j}}&0 \\
    0 & -\left(1-\frac{\eta_{r}+\eta_{1}}{2}\right)\beta & i \omega_{0}+d_{3} & 0 \\
    0 & -\left(\frac{\eta_{r}+\eta_{1}}{2}\right)\beta & 0 & i \omega_{0}+d_{3}\\
  \end{array}
\right)$
$\left(
\begin{array}{c}
    E_{1}^{(1)}\\
    E_{1}^{(2)} \\
    E_{1}^{(3)}  \\
    E_{1}^{(4)} \\
  \end{array}\right) \\=$
 $2\left(
 \begin{array}{c}
 -(1+a)\frac{r}{T_{max}}-(1-c\eta_{1})\alpha b e^{-2 \textrm{i} \omega_{0}\tau_{j}} \\
 (1-c\eta_{1})\alpha b e^{-2 \textrm{i} \omega_{0}\tau_{j}}  \\
 0 \\
 0 \\
 \end{array}
 \right)$\\\\
From above, we can easily calculate a constant vector $E_{1}=(E_{1}^{(1)}, E_{1}^{(2)}, E_{1}^{(3)}, E_{1}^{(4)})\in R^{4}$. Similarly, substituting (\ref{eqn-s21-hep}) and (\ref{eqn-s25-hep}) into (\ref{eqn-s23-hep}), we can get\\\\
$\left(
  \begin{array}{cccc}
   F& -\frac{r T^{*}}{T_{max}}&-(1-c \eta_{1})\alpha T^{*}&0\\
    -(1-c \eta_{1})\alpha V_{I}^{*} &-d_{2} & -(1-c \eta_{1})\alpha T^{*}&0 \\
    0 & \left(1-\frac{\eta_{r}+\eta_{1}}{2}\right)\beta & -d_{3} & 0 \\
    0 & \left(\frac{\eta_{r}+\eta_{1}}{2}\right)\beta & 0 & -d_{3}\\
  \end{array}
\right)$
$\left(
\begin{array}{c}
    E_{2}^{(1)}\\
    E_{2}^{(2)} \\
    E_{2}^{(3)}  \\
    E_{2}^{(4)} \\
  \end{array}\right)\\=$
$ 2 \left(
\begin{array}{c}
-(1+Re(a))\frac{r}{T_{max}}-(1-c\eta_{1})\alpha Re(b) \\
(1-c\eta_{1})\alpha Re(b)  \\
0 \\
0 \\
\end{array}
\right)$\\\\
In the similar manner, we can calculate the constant vector $E_{2}=(E_{2}^{(1)}, E_{2}^{(2)}, E_{2}^{(3)}, E_{2}^{(4)})\in R^{4}.$

\section*{Acknowledgments}

This study was supported by Initiation Grant (Grant number IITR/SRIC/100518) from the Indian Institute of Technology Roorkee, Roorkee, India.

\end{document}